\documentclass[english]{article}
\usepackage[T1]{fontenc}
\usepackage[latin9]{inputenc}
\usepackage{amsthm}
\usepackage{amsmath}
\usepackage{amssymb}
\usepackage{esint}

\makeatletter
\theoremstyle{plain}
\newtheorem{thm}{\protect\theoremname}
  \theoremstyle{plain}
  \newtheorem{cor}[thm]{\protect\corollaryname}
  \theoremstyle{plain}
  \newtheorem{lem}[thm]{\protect\lemmaname}
  \theoremstyle{plain}
  \newtheorem{prop}[thm]{\protect\propositionname}

\makeatother

\usepackage{babel}
  \providecommand{\corollaryname}{Corollary}
  \providecommand{\lemmaname}{Lemma}
  \providecommand{\propositionname}{Proposition}
\providecommand{\theoremname}{Theorem}

\begin{document}

\title{Gradient estimates for porous medium and fast diffusion equations
by martingale method}

\author{Ying Hu\thanks{
IRMAR, Universit\'e Rennes 1, 35042 Rennes Cedex, France. Email:
ying.hu@univ-rennes1.fr.}, Zhongmin Qian\thanks{Mathematical Institute, University of Oxford, Oxford OX2 6GG, England.
Email: qianz@maths.ox.ac.uk.}
, Zichen Zhang\thanks{Mathematical Institute, University of Oxford, Oxford OX2 6GG, England.
Email: zhangz@maths.ox.ac.uk.}}

\date{}

\maketitle
\tableofcontents{}

\begin{abstract}
In this paper, we establish several local and global gradient estimates for the positive
solution of Porous Medium Equations (PMEs) and Fast Diffusion Equations
(FDEs). Our proof is probabilistic and uses martingale techniques.
\end{abstract}

\section{Introduction}

The present paper is devoted to a study of gradient estimates for
solutions of partial differential equations (PDEs) of the form
\begin{equation}
\frac{\partial u}{\partial t}=\Delta u^{m}\label{eq:20130306-1}
\end{equation}
where $m\in\left(0,\infty\right)$ is a given exponent. The problem
is posed on $\left(0,\infty\right)\times\mathbb{R}^{n}$, and  $\Delta$
is the Laplace operator with respect to space variables. When $m=1$,
it is the celebrated heat equation, which is linear and parabolic,
hence enjoying many nice properties. When $m\neq1$, the story becomes
quite different. Let us focus on non-negative solutions to the equation
(\ref{eq:20130306-1}) with initial condition
\[
u\left(0,x\right)=u_{0}\left(x\right),
\]
 the so-called Cauchy problem related to (\ref{eq:20130306-1}), and
collect some results and discuss some phenomena that do not appear
in the case of the  heat equation.

There are two critical values of the exponent $m$, namely $m=1$
and $m=\frac{n-2}{n}$. By rewriting (\ref{eq:20130306-1}) as 
\begin{equation}
\frac{\partial u}{\partial t}=mu^{m-1}\Delta u+m\left(m-1\right)u^{m-2}\left|\nabla u\right|^{2}\label{eq:20130329-1}
\end{equation}
one can see that the coefficient matrix of the second order derivative
is $mu^{m-1}I_{n\times n}$. If $m<1$, this matrix takes the value
of infinity where $u=0$, which means the parabolicity is singular.
While if $m>1$, the matrix vanishes when $u=0$, which means the
parabolicity is degenerate. If we interpret (\ref{eq:20130306-1})
as a differential equation describing a diffusion, this means when
$m<1$, the diffusion is very fast at places where $u$ is small.
Therefore in this case we call (\ref{eq:20130306-1}) the fast diffusion
equation (FDE). If $m>1$, the diffusion slows down wherever $u$
is small, and in this case we call it the porous medium equation (PME).
Notice that since we only consider non-negative solution $u$, from
(\ref{eq:20130329-1}) it is clear that (\ref{eq:20130306-1}) is
always formally parabolic. 

Due to the different behaviors of PME and FDE at $u=0$, theories
about existence and uniqueness of Cauchy problems for these two types
of equations have been studied separately. The commonly used framework
for PME is the $L^{1}\left(\mathbb{R}^{n}\right)$ space. By Theorem
9.12 and Proposition 9.13 in \cite{MR2286292}, if $u_{0}\in L^{1}\left(\mathbb{R}^{n}\right)$,
then there is a unique strong solution which is continuous on $\left(0,\infty\right)\times\mathbb{R}^{n}$.
Moreover, if $u_{0}\in L^{1}\left(\mathbb{R}^{n}\right)$ which is
strictly positive and continuous, then the solution must be smooth.
If we move beyond the scope of $L^{1}\left(\mathbb{R}^{n}\right)$
setting and impose a weaker growth condition on $u_{0}$ such as 
\begin{equation}
\sup_{R\geq1}R^{-\left(n+\frac{2}{m-1}\right)}\int_{\left|x\right|\leq R}\left|u_{0}\left(x\right)\right|dx<\infty,\label{eq:20130903-1}
\end{equation}
then by \cite{MR726106} a unique solution in distribution sense exists
on $\left(0,T\left(u_{0}\right)\right)\times\mathbb{R}^{n}$, where
$T\left(u_{0}\right)\in\left(0,\infty\right]$ depending on $u_{0}$.
By Theorem 3.1 in \cite{MR712265}, the previous growth condition
is satisfied by any continuous non-negative solution. Therefore, condition
(\ref{eq:20130903-1}) is optimal for the class of continuous non-negative
solutions. If the initial data is allowed to be measure-valued, \cite{MR908054}
gives a result which requires similar growth condition as (\ref{eq:20130903-1}).
As for FDE, no requirement on the growth of initial data is needed.
In fact, by Theorem 2.1 and 2.3 in \cite{MR797051}, there exists
a unique solution $u\in C\left(\left[0,\infty\right);L_{loc}^{1}\left(\mathbb{R}^{n}\right)\right)$
in distribution sense if $u_{0}\in L_{loc}^{1}\left(\mathbb{R}^{n}\right)$.
Alternatively, if we impose some growth and decay conditions on $u_{0}$,
by Theorem 1 in \cite{MR1804858} there will be a classical solution
in $\left[0,\bar{T}\right]\times\mathbb{R}^{n}$, and $\bar{T}$ is
finite. 

The degeneracy of parabolicity of PME leads to finite propagation
of its solution, which is one of the special feature of PME. In particular,
by Theorem 14.6 in \cite{MR2286292}, if $u\left(t_{1},\cdot\right)$
is compactly supported in $\mathbb{R}^{n}$, so is $u\left(t_{2},\cdot\right)$
for any $t_{2}>t_{1}$. Consequently, for this kind of solutions,
there is a set in $\left(0,\infty\right)\times\mathbb{R}^{n}$ that
separates the region on which $u$ is positive and the region where
$u$ is zero. According to Theorem 3.3 in \cite{MR570687}, this set,
or the so-called free boundary, is locally H\"older continuous on $\left(0,\infty\right)\times\mathbb{R}^{n}$.
Moreover, in viewing it as a family of boundaries in $\mathbb{R}^{n}$
indexed by $t\in\left(0,\infty\right)$, those boundaries expands
to infinity as $t\rightarrow\infty$ \cite{MR570687}. When the solution
overflows the support of the initial data at some finite time $t_{\star}$,
we see an improvement in the regularity of the free boundary. More
precisely, Theorem 3 in \cite{MR891781} asserts that the free boundary
is locally Lipschitz continuous on $\left(t_{\star},\infty\right)\times\mathbb{R}^{n}$.

Although $m=1$ is a crucial value when talking about finite propagation
and existence theories of equation (\ref{eq:20130306-1}), it is not
a significant value in the study of extinction in finite time and
smoothing effect, where the value $m=\frac{n-2}{n}$ becomes decisive.
The extinction in finite time of a solution is the phenomenon that
arises only when $m<\frac{n-2}{n}$. On page 174 in \cite{MR604277},
it is proved that any solution with initial value $u_{0}\in L^{p_{\star}}\left(\mathbb{R}^{n}\right)\cap L^{1}\left(\mathbb{R}^{n}\right)$,
where $p_{\star}=\frac{n\left(1-m\right)}{2}$, becomes identically
zero after a finite time. More generally, by Theorem 5.2 in \cite{Vazquez2005},
the same result holds for $u_{0}\in{\cal M}^{p_{\star}}$, where 
\[
{\cal M}^{p_{\star}}=\left\{ f\in L_{loc}^{1}\left(\mathbb{R}^{n}\right):\int_{K}\left|f\left(x\right)\right|dx\leq C\left|K\right|^{1-\frac{1}{p_{\star}}},\forall K\ \mbox{with }\left|K\right|<\infty\right\} .
\]
According to Lemma 5.6 in \cite{Vazquez2005}, this is already very
close to the sufficient condition for a solution to extinct in finite
time. One can see from these results that even a positive initial
data may produce a solution that vanishes completely in finite time,
which is quite striking. The reason behind, is the failure of conservation
of mass when $m<\frac{n-2}{n}$, as explained in Section 5.5 in \cite{Vazquez2005}. 

Next, let us describe regularities of solutions in terms of boundedness,
positivity and smoothness. In general, for any $m\in\left(0,\infty\right)$
and $p\in\left[1,\infty\right]$, the solution decreases in $L^{p}\left(\mathbb{R}^{n}\right)$
norm as it evolves in time, according to Theorem 7.2 in \cite{MR2117623}.
Moreover, if $m>\frac{n-2}{n}$, then by Section 3.4 in \cite{Vazquez2005},
for any $p$ and $q$ such that $1\leq p\leq q\leq\infty$,
\begin{equation}
\left\Vert u\left(t,\cdot\right)\right\Vert _{q}\leq c\left(m,n,p,q\right)\left\Vert u_{0}\right\Vert _{p}^{\sigma}t^{-\alpha}\label{eq:20130311-1}
\end{equation}
where 
\[
\alpha=\frac{n\left(q-p\right)}{q\left(n\left(m-1\right)+2p\right)},\ \sigma=\frac{p\left(n\left(m-1\right)+2q\right)}{q\left(n\left(m-1\right)+2p\right)}.
\]
In particular, from this result we see that initial data in $L^{1}\left(\mathbb{R}^{n}\right)$
produce solutions $u\left(t,\cdot\right)\in L^{\infty}\left(\mathbb{R}^{n}\right)$
for any time $t>0$, which is termed as smoothing effect. When $m=\frac{n-2}{n}$,
this is no longer true. Appendix A.3 in \cite{MR1929929} constructed
a solution which is not bounded at any time while still having an
initial value in $L^{1}\left(\mathbb{R}^{n}\right)$. Moreover, when
$m<\frac{n-2}{n}$ and $n\geq3$, Theorem 5.14 in \cite{Vazquez2005}
shows that (\ref{eq:20130311-1}) holds for $q=1$ and $1<p<p_{\star}$.
This means that $L^{p}\left(\mathbb{R}^{n}\right)$ data yield solutions
only in $L^{1}\left(\mathbb{R}^{n}\right)$, a somehow backward smoothing
effect. Nevertheless, we are still able to get a bounded solution
in the case of $m<\frac{n-2}{n}$ if the initial data belongs to a
better space. By Theorem 6.7 in \cite{Vazquez2005}, if $m<\frac{n-2}{n}$,
$p\geq p_{\star}$ and $u_{0}\in L^{p}\left(\mathbb{R}^{n}\right)+L^{\infty}\left(\mathbb{R}^{n}\right)$,
then $u\left(t,\cdot\right)$ is locally bounded and smooth for any
$t>0$. Quantitatively, Theorem 2.1 in \cite{MR2565541} gives a local
upper bound for $u\left(t,\cdot\right)$ in terms of the $L^{p}\left(\mathbb{R}^{n}\right)$
($p>p_{\star}$) norm of $u_{0}$ over a larger local region. 

As for positivity and smoothness, Theorem 3.1 in \cite{MR712265},
Theorem 1.1 in \cite{MR2261689} and Theorem 1.1 in \cite{MR2565541}
give local lower bounds for solutions to (\ref{eq:20130306-1}) in
terms of local $L^{1}\left(\mathbb{R}^{n}\right)$ norm of initial
data in the cases of $m>1$, $m\in\left(\frac{n-2}{n},1\right)$ and
$m\in\left(0,\frac{n-2}{n}\right)$ respectively. These lead to results
about the positivity of solutions. For example, Proposition 1.1 in
\cite{MR891781} gave a necessary and sufficient condition for the
positivity of $u\left(t,x\right)$ when $m>1$ which read as 
\[
\sup_{R>0}R^{-\left(n+\frac{2}{m-1}\right)}\int_{\left|y-x\right|\leq R}\left|u_{0}\left(y\right)\right|dy=\infty.
\]
Besides positivity, if the solution is also locally bounded, then
standard quasilinear theory \cite{MR0241821} implies the smoothness
of the solution on that region. In particular, when $m\in\left(\frac{n-2}{n},1\right)$,
non-negative locally integrable initial data always produce positive
and smooth global solutions, according to the remark after Theorem
2.3 in \cite{MR797051}. This is not true when $m\in\left(0,\frac{n-2}{n}\right)$,
as non-negative locally integrable initial data is not enough for
local boundedness of solutions, which can be seen from solution (0.2)
in \cite{MR2565541}. Moreover, when $m$ is in this range, extinction
in finite time may occur, which kills positivity of solutions in a
global scale. As for PME, in general, solutions are only locally H\"older
continuous. Theorem 4.1 in \cite{MR570687} states that when $u_{0}$
is non-negative, bounded and belongs to $L^{2}\left(\mathbb{R}^{n}\right)$,
$u^{m}$ is uniformly H\"older continuous in every set $\left(\eta_{0},\infty\right)\times\mathbb{R}^{n}$,
$\eta_{0}>0$.  Theorem 7.17 in \cite{MR2286292} tells us that $u$
is locally H\"older continuous on the region where $u$ is bounded. 

During the development of the above works, one of the main tools is
the comparison principle for equation (\ref{eq:20130306-1}), which
is established in Theorem 7.3 in \cite{MR2117623}. Generally speaking,
the comparison is in terms of mass concentration of radially symmetric
functions. Once this comparison of mass concentration is obtained,
comparison in terms of $L^{p}\left(\mathbb{R}^{n}\right)$ ($p\in\left[1,\infty\right]$)
norm follows. There are several special explicit solutions to (\ref{eq:20130306-1})
that are often used in combination with the comparison principle.
Let us only mention three of them here. The first one is source-type
solution, also known as Zel'dovich-Kompanyeets-Barenblatt (ZKB) solution,
which takes the form 
\[
\mathcal{U}\left(t,x\right)=t^{-\alpha}\left(c-k\left|x\right|^{2}t^{-2\beta}\right)_{+}^{\frac{1}{m-1}},
\]
with 
\[
\alpha=\frac{n}{n\left(m-1\right)+2},\ \beta=\frac{\alpha}{n},\ k=\frac{\alpha\left(m-1\right)}{2mn}
\]
and a positive constant $c$. For $m>\frac{n-2}{n}$, 
\[
\lim_{t\rightarrow0}\mathcal{U}\left(t,x\right)=C\left(c,n,m\right)\delta_{0}\left(x\right)
\]
in the sense of measures. So it is the solution to (\ref{eq:20130306-1})
with Dirac delta as initial trace. One can see that it reproduces
the heat kernel as $m\rightarrow1$. Since ZKB solution has compact
support when $m>1$, it plays an important role in the study of finite
propagation speed. When $m\leq\frac{n-2}{n}$, $\mathcal{U}\left(t,x\right)$
is no longer integrable in space variables. Therefore, although still
solving (\ref{eq:20130306-1}), it is not a solution to any Cauchy
problem related to (\ref{eq:20130306-1}). The second type of solution
is a family of self similar solutions, which has the form
\[
\mathcal{U}\left(t,x\right)=t^{-\lambda_{1}}F\left(\left|x\right|t^{-\lambda_{2}}\right).
\]
By a scaling argument, it is shown in section 3.2.1 of \cite{MR2282669}
that any solution to (\ref{eq:20130306-1}) on $\left(0,\infty\right)\times\mathbb{R}^{n}$
with initial data $\left|x\right|^{-\lambda_{3}}$ must have the above
form. Notice that in fact the ZKB solution also belongs to this type.
The third type of solutions we would like to mention is a variance
of the self similar solutions. It has the form
\[
\mathcal{U}\left(t,x\right)=\left(T-t\right)^{-\lambda_{1}}\bar{F}\left(\left|x\right|\left(T-t\right)^{-\lambda_{2}}\right)
\]
with $T>0$. Usually it satisfies (\ref{eq:20130306-1}) only for
$t<T$, after which it blows up or vanishes. For this type of solutions,
probably the most popular and explicit one is 
\[
U\left(x,t;T\right)=2\left(n-\frac{2}{1-m}\right)\left(\frac{T-t}{\left|x\right|^{2}}\right)^{\frac{1}{1-m}}
\]
with $m<\frac{n-2}{n}$, which is a good example of solutions that
extinct in finite time. 

As frequently seen in the field of PDE, a crucial step in the study
of equation (\ref{eq:20130306-1}) is to derive various types of estimates
for solutions. In \cite{MR524760}, Aronson and B\'enilan established
the following gradient estimate for solutions to (\ref{eq:20130306-1}).
If $m\in\left(\frac{n-2}{n},\infty\right)$, $u$ is a positive smooth
solution to (\ref{eq:20130306-1}) and $v=\frac{m}{m-1}u^{m-1}$,
then 
\[
\Delta v\geq-\frac{\alpha}{\left(m-1\right)t}
\]
 where $\alpha=\frac{n\left(m-1\right)}{n\left(m-1\right)+2}$, which
is equivalent to 
\[
\frac{\left|\nabla v\right|^{2}}{\left(m-1\right)v}-\frac{v_{t}}{\left(m-1\right)v}\leq\frac{\alpha}{\left(m-1\right)t}.
\]
This fundamental estimate is then employed in \cite{MR726106} for
the study of existence theory, in \cite{MR797051} for $L_{loc}^{\infty}\left(\mathbb{R}^{n}\right)$
estimate for solutions, and in \cite{MR891781} for obtaining regularity
results for the free boundary of solutions, to name but a few. Later
in \cite{MR2487898}, a local version of Aronson-B\'enilan estimate
was obtained by P. Lu, L. Ni, J. L. Vazquez and C. Villani. They studied
the same problem posed on a local ball of a Riemannian manifold. Let
$B\left(\mathcal{O},2R\right)$ denote a ball with center $\mathcal{O}$
and radius $2R>0$ . Assume that $u$ is a positive solution to (\ref{eq:20130306-1})
on $B\left(\mathcal{O},2R\right)\times\left[0,T\right]$ and the Ricci
curvature Ric$\geq-\left(n-1\right)K^{2}$ on $B\left(\mathcal{O},2R\right)$
for some $K\geq0$. They showed that for any $m>1$ and $\beta>1$,
it holds on $B\left(\mathcal{O},R\right)\times\left[0,T\right]$ that
\[
\frac{\left|\nabla v\right|^{2}}{v}-\beta\frac{v_{t}}{v}\leq\alpha\beta^{2}\left(\frac{1}{t}+C_{2}K^{2}v_{\max}^{2R,T}\right)+\alpha\beta^{2}\frac{v_{\max}^{2R,T}}{R^{2}}C_{1}
\]
where 
\[
v_{\max}^{2R,T}=\max_{B\left(\mathcal{O},2R\right)\times\left[0,T\right]}v.
\]
For $m\in\left(1-\frac{2}{n},1\right)$, they proved that on $B\left(\mathcal{O},R\right)\times\left[0,T\right]$,
for any $\gamma\in\left(0,1\right)$, 
\[
\frac{\left|\nabla v\right|^{2}}{v}-\gamma\frac{v_{t}}{v}\geq\frac{\alpha\gamma^{2}}{C_{3}}\left(\frac{1}{t}+C_{4}\sqrt{C_{3}}K^{2}\bar{v}_{\max}^{2R,T}\right)+\frac{\alpha\gamma^{2}}{C_{3}}\frac{\bar{v}_{\max}^{2R,T}}{R^{2}}C_{5}
\]
where 
\[
\bar{v}_{\max}^{2R,T}=\max_{B\left(\mathcal{O},2R\right)\times\left[0,T\right]}\left(-v\right).
\]
Later in \cite{huang2011gradient} several results of similar type
were obtained by G. Huang, Z. Huang and H. Li. Note that these gradient
bounds do not depend on the initial data. While in \cite{MR1305712},
S.T. Yau established a similar type of gradient bounds depending on
derivatives of initial data for degenerate parabolic equations of
the form
\[
\frac{\partial u}{\partial t}=\Delta\left(F\left(u\right)\right)
\]
with $F\in C^{2}\left(0,\infty\right)$ and $F^{\prime}>0$. In particular,
as explained in \cite{MR2392508}, Yau's result implies that for any
function $c\left(t\right)\in C^{1}\left(0,\infty\right)$ satisfying
\[
\begin{cases}
c\left(t\right)\leq0\\
c^{\prime}\left(t\right)\geq0\\
\left|\nabla v\right|^{2}-2v_{t}+2m\left(\frac{m-1}{m}v\right)^{\frac{m-2}{m-1}}\leq c\left(t\right) & \mbox{at }t=0
\end{cases}
\]
it holds for all $t>0$ that
\[
\left|\nabla v\right|^{2}-2v_{t}+2m\left(\frac{m-1}{m}v\right)^{\frac{m-2}{m-1}}\leq c\left(t\right).
\]

Besides gradient estimates of Aronson-B\'enilan type, the Hamilton type
estimate also plays an important role. It originates from Hamilton
\cite{MR1230276} where it was proved that a positive smooth solution
$u$ to the heat equation on a compact manifold without boundary and
with Ric$\geq-k$, $k>0$, we have
\begin{equation}
\frac{\left|\nabla u\right|^{2}}{u^{2}}\leq\left(\frac{1}{t}+2k\right)\ln\frac{\left\Vert u\right\Vert _{\infty}}{u}.\label{eq:20130605-2}
\end{equation}
This is an upper bound on the gradient of space variables only, hence
leading to a different type of Harnack inequalities. As proved by
Kotschwar \cite{MR2317980}, the same result holds for complete noncompact
manifolds as well. Later a local version was obtained by Souplet and
Zhang in \cite{MR2285258}. As for PME and FDE, L. A. Caffarelli,
J. L. V\'azquez and N. I. Wolanski \cite{MR891781} discussed the case
where the initial data is compactly supported. Namely, under the assumption
that the initial data $u_{0}\geq0$ is integrable and compactly supported,
they established that for $m>1$, there exists a time $T=T\left(u_{0}\right)>0$
and a constant $c=c\left(m,n\right)>0$ such that 
\[
\left|\nabla v\left(x,t\right)\right|\leq c\left(\left(\frac{v}{t}\right)^{\frac{1}{2}}+\frac{\left|x\right|}{t}\right)
\]
for any $t>T$ and almost every $x\in\mathbb{R}^{n}$. Later X. Xu
\cite{MR2853544} derived a local result on a complete Riemannian
manifold with Ric$\geq-k$ for some $k\geq0$. For $m>1$, if there
exists a constant $\delta\in\left(0,\frac{4}{n-1}\right]$ such that
\[
1\leq\frac{v_{\max}^{2R,T}}{v_{\min}^{2R,T}}<\frac{1}{1+\delta}\left(\frac{4m}{\left(n-1\right)\left(m-1\right)}+1\right)
\]
then on $B\left(x_{0},R\right)\times\left[t_{0}-\frac{T}{2},t_{0}\right]$
\[
\frac{\left|\nabla v\right|}{v_{\max}^{2R,T}\left(1+\delta\right)-v}\leq C_{6}\left(m,n\right)\left(\frac{1+\delta}{2\rho\delta R}+\frac{1}{\sqrt{\frac{m-1}{m}v_{\max}^{2R,T}\delta\rho T}}+\sqrt{\frac{k}{\delta}}\right)
\]
where
\[
\rho=2m-\frac{\left(n-1\right)\left(m-1\right)}{2}\frac{v_{\max}^{2R,T}\left(1+\delta\right)-v_{\min}^{2R,T}}{v_{\min}^{2R,T}},
\]
 
\[
v_{\max}^{2R,T}=\sup_{B\left(x_{0},2R\right)\times\left[t_{0}-T,t_{0}\right]}v,\ v_{\min}^{2R,T}=\inf_{B\left(x_{0},2R\right)\times\left[t_{0}-T,t_{0}\right]}v.
\]
 For $m\in\left(1-\frac{4}{n+3},1\right)$, they obtained that on
$B\left(x_{0},R\right)\times\left[t_{0}-\frac{T}{2},t_{0}\right]$,
\[
\frac{\left|\nabla v\right|}{-v}\leq C_{7}\left(m,n\right)\left(\frac{1}{2R}+\frac{1}{\sqrt{\frac{1-m}{m}\bar{v}_{\min}^{2R,T}T}}+\sqrt{k}\right)
\]
where 
\[
\bar{v}_{\min}^{2R,T}=\inf_{B\left(x_{0},2R\right)\times\left[t_{0}-T,t_{0}\right]}\left(-v\right).
\]
This is a generalization of Li Ma et al. \cite{MR2392508}, where
the same estimate was derived only for $n=2$ or $3$ with $m\in\left(1-\frac{1}{\sqrt{n}},1\right)$.
In X. Zhu \cite{Zhu2013201}, it was proved that for $m\in\left(1,1+\frac{1}{\sqrt{2n}+1}\right)$,
on $B\left(x_{0},R\right)\times\left[t_{0}-\frac{T}{2},t_{0}\right]$
\[
v^{\frac{1}{4}\frac{2-m}{m-1}}\left|\nabla v\right|\leq C_{8}\left(v_{\max}^{2R,T}\right)^{1+\frac{1}{4}\frac{2-m}{m-1}}\left(\frac{1}{2R}+\frac{1}{\sqrt{T}}+\sqrt{k}\right)
\]
with
\[
v_{\max}^{2R,T}=\sup_{B\left(x_{0},2R\right)\times\left[t_{0}-T,t_{0}\right]}v.
\]
In X. Zhu \cite{MR2763753}, a gradient bound for $m\in\left(1-\frac{2}{n},1\right)$
was obtained. On $B\left(x_{0},R\right)\times\left[t_{0}-\frac{T}{2},t_{0}\right]$
\[
\frac{\left|\nabla v\right|}{\sqrt{-v}}\leq C_{9}\sqrt{\bar{v}_{\max}^{2R,T}}\left(\frac{1}{2R}+\frac{1}{\sqrt{T}}+\sqrt{k}\right)
\]
where 
\[
\bar{v}_{\max}^{2R,T}=\sup_{B\left(x_{0},2R\right)\times\left[t_{0}-T,t_{0}\right]}\left(-v\right).
\]

In this paper, let us denote by $\mathcal{A}_{x_{0}}^{R,\epsilon}$
the annulus $B\left(x_{0},R+\epsilon\right)\setminus B\left(x_{0},R\right)$.
Set
\[
L=\frac{8}{\epsilon^{2}}\left(\frac{n\epsilon}{2R+\epsilon}+\frac{8\left(R+\epsilon\right)^{2}}{\left(2R+\epsilon\right)^{2}}\right),
\]
 
\begin{equation}
h_{-}=m-\sqrt{1-\left(n-1\right)\left(m-1\right)^{2}},\ h_{+}=m+\sqrt{1-\left(n-1\right)\left(m-1\right)^{2}}.\label{eq:20131212-1}
\end{equation}
We derived several local and global gradient estimates as follows.
\begin{thm}
If $u$ is a positive and bounded solution to \textup{(\ref{eq:20130306-1})}
on \textup{$\left(0,t_{1}\right)\times B\left(x_{0},R+\epsilon\right)$}
with $m\in\left(1-\frac{1}{\sqrt{n-1}},1+\frac{1}{\sqrt{n-1}}\right)$,
then for any $h\in\left(h_{-},h_{+}\right),$ $\rho\in\left(0,m\left|\frac{2m+h}{2h}\right|^{-1}\right)$,
and $\left(T,x\right)\in\left(0,t_{1}\right)\times B\left(x_{0},R\right)$,
we have \textup{
\begin{eqnarray}
 &  & u^{m+h-3}\left|\nabla u\right|^{2}\left(T,x\right)\nonumber \\
 & \leq & \frac{\left\Vert I_{B\left(x_{0},R+\epsilon\right)}u_{0}^{h}\right\Vert _{\infty}-u^{h}\left(T,x\right)}{hm\left(m-\rho\left|\frac{2m+h}{2h}\right|\right)T}+\frac{\left(\rho+\left|\frac{2m+h}{2}\right|\right)LM_{2}}{\left|h\right|\rho\left(m-\rho\left|\frac{2m+h}{2h}\right|\right)}\nonumber \\
 &  & +\left(\left|2m-h\right|l^{-3}+6l^{-1}\right)mM_{1}L^{2}\frac{T}{8}\label{eq:20130528-1}
\end{eqnarray}
}where \textup{
\[
l=\frac{\left(h-h_{-}\right)\left(h_{+}-h\right)}{3\left|m-\frac{h}{2}\right|+3},\ M_{1}=\sup_{\left[0,T\right]\times\mathcal{A}_{x_{0}}^{R,\epsilon}}u^{2m+h-2},\ M_{2}=\sup_{\left[0,T\right]\times\mathcal{A}_{x_{0}}^{R,\epsilon}}u^{m+h-1}.
\]
}
\end{thm}
Comparing to the existing work, our result is essentially about estimating
the gradient at time $t_{0}$ with the information of $u$ during
$\left[t_{0}-T,t_{0}\right]$, while the other results we listed above
are using information of $u$ during $\left[t_{0}-2T,t_{0}\right]$.
Therefore, we have an extra term linear in $T$. In terms of the constraint
on $m$, one can see that when the dimension $n=1$, there is no constraint.
When $n>1$, comparing to all the existing results we are aware of,
our bound is valid for a wider range of $m$ when $m<1$. As for the
case when $m>1$, our range of $m$ is larger than that in \cite{Zhu2013201}.
If we look for a bound on the pressure variable, by taking $h=0$,
the restriction on $m$ becomes identical to the one in \cite{MR2763753}.
\begin{cor}
If $u$ is a positive and bounded solution to \textup{(\ref{eq:20130306-1})}
on \textup{$\left(0,t_{1}\right)\times B\left(x_{0},R+\epsilon\right)$}
with $m\in\left(1-\frac{2}{n},1\right)$, and the pressure variable
$v=\frac{u^{m-1}}{m-1}$, then for any $\rho\in\left(0,m\right)$,
and $\left(T,x\right)\in\left(0,t_{1}\right)\times B\left(x_{0},R\right)$,
we have
\begin{eqnarray*}
\frac{\left|\nabla v\right|^{2}}{\left(m-1\right)v}\left(T,x\right) & \leq & \frac{\left\Vert I_{B\left(x_{0},R+\epsilon\right)}\log u_{0}\right\Vert _{\infty}-\log u\left(T,x\right)}{m\left(m-\rho\right)T}\\
 &  & +\frac{LM_{3}}{\rho\left(m-\rho\right)}+\left(\frac{m}{2}l^{-3}+\frac{3}{2}l^{-1}\right)mM_{1}L^{2}\frac{T}{2}
\end{eqnarray*}
 where 
\[
l=\frac{-n\left(m-1\right)^{2}-2\left(m-1\right)}{3\left(m+1\right)},\ M_{1}=\sup_{\left[0,T\right]\times\mathcal{A}_{x_{0}}^{R,\epsilon}}u^{2m-2},
\]
\textup{
\[
M_{3}=\sup_{\left[0,T\right]\times\mathcal{A}_{x_{0}}^{R,\epsilon}}\left(\rho\log u+\left(m\log u+1\right)^{2}\right)u^{m-1}
\]
}
\end{cor}
As for the global case, by letting $\epsilon=R$, then $R\rightarrow\infty$
and $\rho\rightarrow0$ in (\ref{eq:20130528-1}), we obtain a very
neat gradient bound. 
\begin{cor}
If $u$ is a positive and bounded solution to (\ref{eq:20130306-1})
on \textup{$\left(0,\infty\right)\times\mathbb{R}^{n}$} and $m\in\left(1-\frac{1}{\sqrt{n-1}},1+\frac{1}{\sqrt{n-1}}\right)$,\textup{
}then for any $h\in\left[h_{-},h_{+}\right]$ and $\left(T,x\right)\in\left(0,\infty\right)\times\mathbb{R}^{n}$,
it holds that \textup{
\begin{equation}
u^{m+h-3}\left|\nabla u\right|^{2}\left(T,x\right)\leq\frac{\left\Vert u_{0}^{h}\right\Vert _{\infty}-u^{h}\left(T,x\right)}{Thm^{2}}.\label{eq:20130606-1-1}
\end{equation}
}
\end{cor}
First of all, this bound does not depend on the dimension $n$ explicitly.
This is because the information about dimension has been incorporated
into the initial value and the constraint on $h$. In fact, \cite{MR1044287}
derived a gradient bound which is independent of the initial data,
but dependent on dimension $n$. Secondly, it naturally shows that
when a function touches its maximum, its gradient vanishes. In this
sense it is a tight gradient bound. Moreover, when $m=1$
and $h\rightarrow 0^+$, our results recovers Hamilton's gradient bound
for heat equation. 
\begin{cor}
(Hamilton's gradient estimate) If $u$ is a solution to the heat equation
on \textup{$\left(0,\infty\right)\times\mathbb{R}^{n}$} with positive
and bounded initial data, then for any $\left(T,x\right)\in\left(0,\infty\right)\times\mathbb{R}^{n}$,
it holds that \textup{
\begin{equation}
\frac{\left|\nabla u\right|^{2}}{u^{2}}\left(T,x\right)\leq\frac{1}{T}\log\frac{\left\Vert u_{0}\right\Vert _{\infty}}{u\left(T,x\right)}.\label{eq:20130606-2-1}
\end{equation}
Our last result is an extension of Corollary 3, where the condition
that 
\[
\left(h-h_{-}\right)\left(h-h_{+}\right)<0,
\]
is essential. }

\textup{We next consider the case where $m$ and $h$ satisfy the
condition that 
\[
\left(h-h_{-}\right)\left(h-h_{+}\right)>0.
\]
The previous inequality is equivalent to that $m\in\left(0,1-\frac{1}{\sqrt{n-1}}\right)\cup\left(1+\frac{1}{\sqrt{n-1}},\infty\right)$
and $h\in\mathbb{R}$ , or $m\in\left[1-\frac{1}{\sqrt{n-1}},1+\frac{1}{\sqrt{n-1}}\right]$
but $h\in\left(-\infty,h_{-}\right)\cup\left(h_{+},\infty\right)$.
For such $m$ and $h$, let us define
\begin{equation}
U_{\pm}=\frac{mh\pm\left|h\right|\sqrt{m^{2}+\left(h-h_{-}\right)\left(h-h_{+}\right)}}{\left(h-h_{-}\right)\left(h-h_{+}\right)}.\label{eq:20131212-2}
\end{equation}
}\end{cor}
\begin{thm}
Let $u$ be a positive and bounded solution to (\ref{eq:20130306-1})
on \textup{$\left(0,\infty\right)\times\mathbb{R}^{n}$.} Suppose
$h\in\mathbb{R}$ satisfying the condition that 
\[
\left(h-h_{-}\right)\left(h-h_{+}\right)>0,
\]
and 
\[
U\equiv\frac{\left\Vert u_{0}^{h}\right\Vert _{\infty}}{u_{\min}^{h,T}}-1\in\left(U_{-},U_{+}\right)
\]
where
\[
u_{\min}^{h,T}=\inf_{\left(t,x\right)\in\left[0,T\right]\times\mathbb{R}^{n}}u^{h}\left(t,x\right).
\]
Then
\begin{eqnarray*}
 &  & u^{m+h-3}\left|\nabla u\right|^{2}\left(T,x\right)\\
 & \leq & 2\frac{\left\Vert u_{0}^{h}\right\Vert _{\infty}-u^{h}\left(T,x\right)}{mT\left(h-h_{-}\right)\left(h-h_{+}\right)U^{-1}\left(U-U_{-}\right)\left(U_{+}-U\right)}.
\end{eqnarray*}

\end{thm}
This result, together with Corollary 3 shows that when $m\in\left(1-\frac{1}{\sqrt{n-1}},1+\frac{1}{\sqrt{n-1}}\right)$,
a gradient bound only depending on the maximum of initial data can
be established. While when $m$ is outside of this range, the gradient
bound will depend on the minimum of $u$ as well. A similar phenomena
has been shown in \cite{MR2853544} , but only for $m>1$. Nevertheless,
the meaning behind this phenomena is not clear to us yet, which is
worthy of exploring. 

We would like to mention that employing martingale theory to derive
gradient estimate is not new, and there is a large number of papers
devoted to the study of solutions of PDEs by using stochastic differential
equations, for example \cite{MR1902192}, \cite{MR2605409} and literature
therein for a small sample. On the other hand, to the best of authors'
knowledge, there are few papers dealing with the kind of nonlinear
PDEs by using martingale methods. 

The rest of the paper is organized as follows. In Section 2, we illustrate our main idea by establishing the gradient estimate for solution to the heat equation. In Section 3,
we prove our local and global estimates.

\section{Gradient estimate for solutions to the heat equation}

To illustrate the main idea, let us consider the heat equation in
this chapter. Assume $u$ solves
\begin{eqnarray}
\frac{\partial u}{\partial t} & = & \Delta u\mbox{ on}\left(0,\infty\right)\times\mathbb{R}^{n}\nonumber \\
u\left(x,0\right) & = & u_{0}\left(x\right)\mbox{ on }\mathbb{R}^{n}.\label{eq:20130318-1}
\end{eqnarray}
To avoid technical difficulties, we also assume $u$ is smooth, bounded
and has bounded derivatives with respect to space variables up to
the second order. It is known that there is a close link between a
large class of parabolic PDEs and diffusion processes, in the sense
that the differential operator for a PDE can be identified as a generator
for a diffusion process. Once this one-to-one correspondence has been
established, we are given a way to study a PDE through its diffusion
process counterpart, or the other way around. For example, the solution
to a PDE can be expressed in terms of the expectation of a diffusion
process at a certain time, the so-called path integration, and the
transition probability density function of a diffusion process is
the kernel of a PDE. For a more comprehensive account of this area,
we refer to the book \cite{MR532498} by Stroock and Varadhan. Here
we also relates our PDE, (\ref{eq:20130318-1}), to a diffusion process
in the way that we have just explained. For a given point $\left(T,x\right)\in\left(0,\infty\right)\times\mathbb{R}^{n}$,
define an $n-$dimensional stochastic process $X_{t}$ by solving
the stochastic differential equation
\begin{eqnarray*}
dX_{t}^{\alpha} & = & \sqrt{2}dW_{t}^{\alpha}\\
X_{0}^{\alpha} & = & x^{\alpha}
\end{eqnarray*}
where $\alpha\in\left\{ 1,\cdots,n\right\} $ and $W=\left(W^{1},\cdots,W^{n}\right)$
is a standard $n-$dimensional Brownian motion on the probability
space $\left(\Omega,\mathcal{F},\mathbb{P}\right)$. Then we have
a progressively measurable function $X$ from $\left(\left(0,\infty\right)\times\Omega,\mathcal{B}\left(\left(0,\infty\right)\right)\otimes\mathcal{F}\right)$
to $\mathbb{R}^{n}$. Bearing in mind that the aim is to get an upper
bound for $\left|\nabla u\right|^{2}\left(T,x\right)$, let us consider
the process $\left|\nabla u\right|^{2}\left(T-t,X_{t}\right)$ with
index $t\in\left[0,T\right]$, that is, the gradient running backward
on diffusion process $X$. The way we composite $\left|\nabla u\right|^{2}$
with $X$ is commonly seen. In terms of computation, this will lead
to the disappearance of terms containing the time derivative when
using Ito formula. Intuitively, this is because on the one hand our
diffusion process always starts at the deterministic point that we
are interested in and then evolves in a certain random way as $t$
increases, but on the other hand, due to the nature of parabolic PDEs,
we need to use the information about the solution before $T$. So
we need the time variable to decrease when $t$ increases. 

One can observe that $\left|\nabla u\right|^{2}\left(T-t,X_{t}\right)$
is a semimartingale. Hence we are encouraged to turn to the theory
of martingales, a concept introduced by Paul L\'evy in 1930s and greatly
developed by J.L Doob in his book \cite{MR0058896}. First of all,
let us decompose the semimartingale $\left(\left|\nabla u\right|^{2}\left(T-t,X_{t}\right)\right)_{0\leq t\leq T}$
into a sum of a local martingale and a process with finite variation.
From now on, when there is no potential confusion, we omit the specification
of variables in functions, as we always consider functions running
backward on diffusion $X$. Taking derivatives with respect to $x^{\alpha}$
on both sides of (\ref{eq:20130318-1}) yields
\begin{equation}
\frac{\partial^{2}u}{\partial t\partial x^{\alpha}}=\Delta\frac{\partial u}{\partial x^{\alpha}}.\label{eq:20130318-2}
\end{equation}
Since $u$ is smooth by our assumption, we can apply Ito formula on
$\frac{\partial u}{\partial x^{\alpha}}\left(T-t,X_{t}\right)$ and
obtain 
\begin{eqnarray*}
d\frac{\partial u}{\partial x^{\alpha}} & = & -\frac{\partial^{2}u}{\partial t\partial x^{\alpha}}dt+\Delta\frac{\partial u}{\partial x^{\alpha}}dt+\sqrt{2}\frac{\partial^{2}u}{\partial x^{\alpha}\partial x^{\beta}}dW_{t}^{\beta}\\
 & = & \sqrt{2}\frac{\partial^{2}u}{\partial x^{\alpha}\partial x^{\beta}}dW_{t}^{\beta}
\end{eqnarray*}
where the last equality results from (\ref{eq:20130318-2}). Then,
again by Ito formula, 
\begin{eqnarray}
d\left|\nabla u\right|^{2} & = & 2\frac{\partial u}{\partial x^{\alpha}}d\frac{\partial u}{\partial x^{\alpha}}+\sum_{\alpha}d\left\langle \frac{\partial u}{\partial x^{\alpha}}\right\rangle \nonumber \\
 & = & 2\sqrt{2}\frac{\partial u}{\partial x^{\alpha}}\frac{\partial^{2}u}{\partial x^{\alpha}\partial x^{\beta}}dW_{t}^{\beta}+2\sum_{\alpha,\beta=1}^{n}\left(\frac{\partial^{2}u}{\partial x^{\alpha}\partial x^{\beta}}\right)^{2}dt.\label{eq:20130320-3}
\end{eqnarray}
This is the decomposition we are looking for. It can be seen that
the finite variation term is non-negative, and the local martingale
part is a true martingale under $\mathbb{P}$, a consequence of the
boundedness of $\frac{\partial u}{\partial x^{\alpha}}$ and $\frac{\partial^{2}u}{\partial x^{\alpha}\partial x^{\beta}}$,
$\forall\alpha,\beta\in\left\{ 1,\cdots,n\right\} $. This means that
$\left|\nabla u\right|^{2}\left(T-t,X_{t}\right)$ is a submartingale,
which implies that
\[
E\left[\left|\nabla u\right|^{2}\left(T-t,X_{t}\right)\right]\geq E\left[\left|\nabla u\right|^{2}\left(T,X_{0}\right)\right]=\left|\nabla u\right|^{2}\left(T,x\right)
\]
 for any $t\in\left[0,T\right]$, where the equality results from
the fact that $X_{0}=x$ $\mathbb{P}-$almost surely. By integrating
both sides on $t$ from $0$ to $T$, we have 
\begin{equation}
\int_{0}^{T}E\left[\left|\nabla u\right|^{2}\left(T-t,X_{t}\right)\right]dt\geq T\left|\nabla u\right|^{2}\left(T,x\right).\label{eq:20130318-7}
\end{equation}
On the other hand, through a scaling argument, one can see that a
bound on $\left|\nabla u\right|^{2}\left(T,x\right)$ involving $T$
and $\left\Vert u\right\Vert _{\infty}$ should have the order $O\left(\frac{\left\Vert u\right\Vert _{\infty}^{2}}{T}\right)$.
Therefore, let us consider the process $u^{2}\left(T-t,X_{t}\right)$.
By Ito formula, 
\begin{equation}
u^{2}\left(0,X_{T}\right)-u^{2}\left(T,X_{0}\right)=2\sqrt{2}\int_{0}^{T}u\frac{\partial u}{\partial x^{\alpha}}dW_{t}^{\alpha}+2\int_{0}^{T}\left|\nabla u\right|^{2}dt.\label{eq:20130318-8}
\end{equation}
Since $u$ and $\frac{\partial u}{\partial x^{\alpha}}$, $\alpha\in\left\{ 1,\cdots,n\right\} $,
are bounded by our assumption, $\int_{0}^{\cdot}u\frac{\partial u}{\partial x^{\alpha}}dW_{t}^{\alpha}$
is a true martingale. Hence $E\left[\int_{0}^{T}u\frac{\partial u}{\partial x^{\alpha}}dW_{t}^{\alpha}\right]=0$.
It then follows from (\ref{eq:20130318-8}) that 
\begin{equation}
E\left[u^{2}\left(0,X_{T}\right)\right]-u^{2}\left(T,x\right)=2E\left[\int_{0}^{T}\left|\nabla u\right|^{2}dt\right].\label{eq:20130318-6}
\end{equation}
By Fubini's theorem, $E\left[\int_{0}^{T}\left|\nabla u\right|^{2}dt\right]=\int_{0}^{T}E\left[\left|\nabla u\right|^{2}\right]dt$.
Therefore, (\ref{eq:20130318-7}) and (\ref{eq:20130318-6}) together
yields
\begin{eqnarray*}
\left|\nabla u\right|^{2}\left(T,x\right) & \leq & \frac{E\left[u^{2}\left(0,X_{T}\right)\right]-u^{2}\left(T,x\right)}{2T}.
\end{eqnarray*}
Since the law of $X_{T}$ is absolutely continuous with respect to
Lebesgue measure, $u^{2}\left(0,X_{T}\right)\leq\left\Vert u_{0}\right\Vert _{\infty}^{2}$
$\mathbb{P}-$ almost surely, where $\left\Vert \cdot\right\Vert _{\infty}$
denotes the essential supremum. Hence we have the following result. 
\begin{thm}
If $u$ solves the Cauchy problem (\ref{eq:20130318-1}) and $u$
is smooth, bounded and has bounded derivatives with respect to space
variables up to the second order, then for any $\left(T,x\right)\in\left(0,\infty\right)\times\mathbb{R}^{n}$,
we have \textup{
\begin{equation}
\left|\nabla u\right|^{2}\left(T,x\right)\leq\frac{\left\Vert u_{0}\right\Vert _{\infty}^{2}-u^{2}\left(T,x\right)}{2T}\leq\frac{\left\Vert u_{0}\right\Vert _{\infty}^{2}}{T}.\label{eq:20130604-1}
\end{equation}
}
\end{thm}
We want to remark here that this method also applies to estimates
of gradients of higher orders. To see this, first observe that by
considering the function $\bar{u}\left(t,x\right)=u\left(t+\epsilon,x\right)$
where $\epsilon>0$ in the above argument, we get 
\[
\left|\nabla u\right|^{2}\left(T,x\right)\leq\frac{\left\Vert u\left(\epsilon,\cdot\right)\right\Vert _{\infty}^{2}}{T-\epsilon}.
\]
Since $\frac{\partial u}{\partial x^{i}}$ also satisfies the heat
equation, it holds that
\[
\left|\nabla\frac{\partial u}{\partial x^{i}}\right|^{2}\left(T,x\right)\leq\frac{\left\Vert \frac{\partial u}{\partial x^{i}}\left(\epsilon,\cdot\right)\right\Vert _{\infty}^{2}}{T-\epsilon}.
\]
Summing over the index $i$ yields
\[
\sum_{i,j=1}^{n}\left|\frac{\partial^{2}u}{\partial x^{i}\partial x^{j}}\right|^{2}\left(T,x\right)\leq\frac{\left\Vert \left|\nabla u\right|^{2}\left(\epsilon,\cdot\right)\right\Vert _{\infty}}{T-\epsilon}\leq\frac{\left\Vert u_{0}\right\Vert _{\infty}^{2}}{\epsilon\left(T-\epsilon\right)},
\]
where the last inequality results from (\ref{eq:20130604-1}). Then
we can minimize the right hand side by choosing $\epsilon=\frac{T}{2}$,
hence obtaining
\[
\sum_{i,j=1}^{n}\left|\frac{\partial^{2}u}{\partial x^{i}\partial x^{j}}\right|^{2}\left(T,x\right)\leq\frac{4\left\Vert u_{0}\right\Vert _{\infty}^{2}}{T^{2}}.
\]

{\bf Remark.}
An approach based on the link between Backward Stochastic Differential Equation (BSDE) and PDE to establish gradient estimates for positive solutions to the heat equation of elliptic or subelliptic operators on Euclidean spaces or on Riemannian manifolds is developed in \cite{HuQian}.

\section{Gradient estimate for positive solutions to PME and FDE}

\subsection{Local gradient estimate}

Now we move to the equation (\ref{eq:20130306-1}). In this section
we study the problem on a local scale. Denote by $B\left(x_{0},R+\epsilon\right)$
a closed ball in $\mathbb{R}^{n}$ with center $x_{0}$ and radius
$R+\epsilon$, where $R$, $\epsilon>0$. Let $u$ be a positive and
bounded solution to (\ref{eq:20130306-1}) on $\left[0,t_{1}\right]\times B\left(x_{0},R+\epsilon\right)$.
The positivity of $u$ ensures that no degeneracy of parabolicity
would happen in $\left[0,t_{1}\right]\times B\left(x_{0},R+\epsilon\right)$.
Hence we can use theory about non-degenerate quasilinear parabolic
PDE to obtain that $u$ is smooth in $\left[0,t_{1}\right]\times B\left(x_{0},R+\epsilon\right)$.
Note that at this moment $u$ is defined only on $\left[0,t_{1}\right]\times B\left(x_{0},R+\epsilon\right)$.
This brings difficulty to our martingale method, as we will consider
$u$ running backward on a stochastic process, which takes values
on the whole space $\mathbb{R}^{n}$ at any time. To get around this
obstacle, let $\tilde{u}$ be a positive and smooth function with
bounded derivatives of all orders defined on $\left[0,t_{1}\right]\times\mathbb{R}^{n}$,
such that 
\[
u=\tilde{u}\mbox{ on }\left[0,t_{1}\right]\times B\left(x_{0},R+\epsilon\right).
\]
Note that as $u$ is strictly positive and smooth on the compact set
$\left[0,t_{1}\right]\times B\left(x_{0},R+\epsilon\right)$, such
$\tilde{u}$ exists. It is worthwhile to point out here that the behavior
of the extended function $\tilde{u}$ outside $\left[0,t_{1}\right]\times B\left(x_{0},R+\epsilon\right)$
will not enter into our computation in the sequel, as Li-Yau's localization
technique will be adopted. Next, we take a transformation on $\tilde{u}$
by setting 
\begin{equation}
f=\frac{\tilde{u}^{h}}{h}\label{eq:20130320-2}
\end{equation}
 for some $h\in\mathbb{R}\setminus\left\{ 0\right\} $. This is a
generalization of the transform $v=\frac{u^{m-1}}{m-1}$, which repeatedly
appears in literature concerning PME and FDE, such as \cite{MR524760},
\cite{MR2487898}, \cite{MR891781}, \cite{MR2853544} and \cite{MR2763753}.
From (\ref{eq:20130306-1}) one can derive that the so-called pressure
variable $v$ satisfies
\begin{equation}
\frac{\partial v}{\partial t}=m\left(m-1\right)v\Delta v+m\left|\nabla v\right|^{2}.\label{eq:20130320-1}
\end{equation}
We can see that the exponent $m$ in (\ref{eq:20130306-1}) comes
down into coefficients in (\ref{eq:20130320-1}), and both terms on
the right hand side of (\ref{eq:20130320-1}) are quadratic, which
facilitates many computations. However, this feature is not crucial
to our method. Hence we attempt to generalize this transform with
(\ref{eq:20130320-2}). It turns out that the flexibility in choosing
$h$ results in an enlargement of the range of $m$ that our gradient
bound is valid for. 

As we are looking for a gradient bound on a local scale, we are keen
to only use the local information about $f$. For this purpose, we
adopt the localization technique of Li and Yau \cite{MR834612} to
introduce a cut-off function $\phi\in C^{2}\left(\mathbb{R}^{n}\right)$
satisfying 
\[
\phi\left(x\right)=\begin{cases}
1 & \mbox{on }B\left(x_{0},R\right)\\
0 & \mbox{on }B\left(x_{0},R+\epsilon\right)^{c},
\end{cases}
\]
 
\begin{equation}
\left|\Delta\phi\right|\leq L\phi^{\frac{1}{2}}\label{eq:2012-11-28-2}
\end{equation}
and
\begin{equation}
\left|\nabla\phi\right|^{2}\leq L\phi^{\frac{3}{2}},\label{eq:2012-11-28-3}
\end{equation}
for some $L>0$. Note that such cut-off function exists. One possible
choice is 
\[
\phi\left(x\right)=\begin{cases}
1 & \mbox{on }B\left(x_{0},R\right)\\
0 & \mbox{on }B\left(x_{0},R+\epsilon\right)^{c}\\
\left(\left(\frac{R+\epsilon}{R}\right)^{2}-1\right)^{-4}\left(\left(\frac{R+\epsilon}{R}\right)^{2}-\left|\frac{x-x_{0}}{R}\right|^{2}\right)^{4} & \mbox{otherwise}
\end{cases}
\]
and 
\[
L=\frac{8}{\epsilon^{2}}\left(\frac{n\epsilon}{2R+\epsilon}+\frac{8\left(R+\epsilon\right)^{2}}{\left(2R+\epsilon\right)^{2}}\right).
\]
Then function $\phi f$, the multiplication of functions $\phi$ and
$f$, cuts all the information of $f$ outside $\left[0,t_{1}\right]\times B\left(x_{0},R+\epsilon\right)$,
while faithfully preserving its behavior in $\left[0,t_{1}\right]\times B\left(x_{0},R\right)$.
Therefore, it is this function that we are going to consider in the
sequel.

\subsubsection{From PDE to SDE}

Just as the case of heat equation, let us begin by fixing a point
$\left(T,x\right)\in\left[0,t_{1}\right]\times B\left(x_{0},R\right)$.
From (\ref{eq:20130306-1}) and the definition of $f$, we have on
$\left[0,t_{1}\right]\times\mathbb{R}^{n}$, 
\begin{eqnarray}
\frac{\partial\phi f}{\partial t} & = & m\left(hf\right)^{\frac{m-1}{h}}\phi\Delta f+m\left(m-h\right)\left(hf\right)^{\frac{m-h-1}{h}}\phi\left|\nabla f\right|^{2}\nonumber \\
 & = & m\left(hf\right)^{\frac{m-1}{h}}\Delta\left(\phi f\right)-m\left(hf\right)^{\frac{m-1}{h}}f\Delta\phi\nonumber \\
 &  & -2m\left(hf\right)^{\frac{m-1}{h}}\nabla f\cdot\nabla\phi+m\left(m-h\right)\left(hf\right)^{\frac{m-h-1}{h}}\phi\left|\nabla f\right|^{2}.\label{eq:20120906-1}
\end{eqnarray}
Then let us link PDE (\ref{eq:20120906-1}) with the diffusion process
$X=\left(X_{t}\right)_{0\leq t\leq T}$, whose generator $\mathcal{L}$
is given by 
\[
\mathcal{L}_{t}w\left(y\right)=m\left(hf\right)^{\frac{m-1}{h}}\left(T-t,y\right)\Delta w\left(y\right),\ \forall w\in C_{0}^{2}\left(\mathbb{R}^{n}\right).
\]
$\mathcal{L}$ corresponds to the principle part of the differential
operator in (\ref{eq:20120906-1}). Note that by the definition of
$f$, $\left(hf\right)^{\frac{m-1}{h}}\geq0$ since $u\geq0$. By
\cite{MR532498} the way to obtain $X$ is to solve the stochastic
differential equation (SDE)
\begin{eqnarray*}
dX_{t}^{\alpha} & = & \sqrt{2m}\left(hf\right)^{\frac{m-1}{2h}}\left(T-t,X_{t}\right)dW_{t}^{\alpha}\\
X_{0} & = & x
\end{eqnarray*}
for $t\in\left[0,T\right]$, where $W=\left(W^{1},\cdots,W^{n}\right)$
is a standard n-dimensional Brownian motion on a probability space
$\left(\Omega,\mathcal{F},\mathbb{P}\right)$ and the stochastic integral
is in Ito's sense. By our assumption on $u$, $\left(hf\right)^{\frac{m-1}{2h}}$
is bounded, smooth and has bounded derivatives. Therefore, by \cite{MR1121940}
this SDE has a unique strong solution on the time interval $\left[0,T\right]$.

\subsubsection{Fundamental decompositions}

The next step is to consider various functions running backward on
process $X$. They are $f\left(T-t,X_{t}\right)$, $\phi\left(X_{t}\right)$,
$\left|\nabla f\right|^{2}\left(T-t,X_{t}\right)$ and so on. In fact,
we will omit the specification of variables in functions, as all the
functions below are compositions on $\left(T-t,X_{t}\right)$.  $t$
will always take values in $\left[0,T\right]$, which is exactly the
time interval where $X$ lives. One can see that these processes are
all semimartingales, waiting for us to decompose and then releasing
information. But unlike the heat equation, (\ref{eq:20120906-1})
is not linear, a consequence of the non-linearity of (\ref{eq:20130306-1}).
Therefore, we should be very careful in choosing semimartingales for
decomposition. To begin with, let us consider 
\[
Y_{t}\triangleq\left(\phi f\right)\left(T-t,X_{t}\right).
\]
By Ito formula,
\begin{eqnarray}
dY_{t} & = & \sqrt{2m}\left(hf\right)^{\frac{m-1}{2h}}\frac{\partial\left(f\phi\right)}{\partial x^{\alpha}}dW_{t}^{\alpha}-m\left(m-h\right)\left(hf\right)^{\frac{m-h-1}{h}}\phi\left|\nabla f\right|^{2}dt\nonumber \\
 &  & +\left(\Delta\phi m\left(hf\right)^{\frac{m-1}{h}}f+2m\left(hf\right)^{\frac{m-1}{h}}\frac{\partial f}{\partial x^{\alpha}}\frac{\partial\phi}{\partial x^{\alpha}}\right)dt.\label{eq:20120906-5}
\end{eqnarray}
It is readily seen that 
\begin{equation}
Y_{T}=\phi f\left(0,X_{T}\right)=\phi\frac{u_{0}^{h}}{h}\left(X_{T}\right).\label{eq:2012-11-29-1}
\end{equation}
Observe that in the finite variation part of (\ref{eq:20120906-5}),
the term with highest degree in $\nabla f$ is 
\[
-m\left(m-h\right)\left(hf\right)^{\frac{m-h-1}{h}}\phi\left|\nabla f\right|^{2}\left(T-t,X_{t}\right)dt.
\]
This means that it can be controlled in terms of $\left\Vert u_{0}\right\Vert _{\infty}$
with the help of (\ref{eq:20120906-5}). Therefore, it is this term
that worth investigation. Define 
\[
H_{t}=\left(hf\right)^{\frac{m-h-1}{h}}\phi\left|\nabla f\right|^{2}\left(T-t,X_{t}\right),\ t\in\left[0,T\right].
\]
 The function $\left(hf\right)^{\frac{m-h-1}{h}}\phi\left|\nabla f\right|^{2}$
is smooth on $\left[0,t_{1}\right]\times\mathbb{R}^{n}$. Its dynamic
only depends on the behavior of $u$ within $\left[0,t_{1}\right]\times B\left(x_{0},R+\epsilon\right)$,
which is governed by (\ref{eq:20130306-1}). Therefore, we can use
Ito formula and decompose $H$ into local martingale part and finite
variation part as follows.

\begin{eqnarray}
dH_{t} & = & d\left(\left(hf\right)^{\frac{m-h-1}{h}}\phi\sum_{\alpha}\left|\frac{\partial f}{\partial x^{\alpha}}\right|^{2}\right)\nonumber \\
 & = & 2\sqrt{2m}\left(hf\right)^{\frac{3m-2h-3}{2h}}\phi\frac{\partial^{2}f}{\partial x^{\alpha}\partial x^{\beta}}\frac{\partial f}{\partial x^{\alpha}}dW_{t}^{\beta}\nonumber \\
 &  & +\sqrt{2m}\left(m-h-1\right)\left(hf\right)^{\frac{3m-4h-3}{2h}}\phi\left|\nabla f\right|^{2}\frac{\partial f}{\partial x^{\beta}}dW_{t}^{\beta}\nonumber \\
 &  & +\sqrt{2m}\left(hf\right)^{\frac{3m-2h-3}{2h}}\frac{\partial\phi}{\partial x^{\beta}}\left|\nabla f\right|^{2}dW_{t}^{\beta}\nonumber \\
 &  & -m\left(2m-h+1\right)\left(m-h-1\right)\left(hf\right)^{\frac{2m-3h-2}{h}}\phi\left|\nabla f\right|^{4}dt\nonumber \\
 &  & +2m\left(m-h-1\right)\left(hf\right)^{\frac{2m-2h-2}{h}}\frac{\partial\phi}{\partial x^{\alpha}}\frac{\partial f}{\partial x^{\alpha}}\left|\nabla f\right|^{2}dt\nonumber \\
 &  & +m\left(hf\right)^{\frac{2m-h-2}{h}}\Delta\phi\left|\nabla f\right|^{2}dt\nonumber \\
 &  & -2m\left(m-1\right)\left(hf\right)^{\frac{2m-2-2h}{h}}\phi\Delta f\left|\nabla f\right|^{2}dt\nonumber \\
 &  & -4m\left(hf\right)^{\frac{2m-2h-2}{h}}\phi\frac{\partial f}{\partial x^{\beta}}\frac{\partial^{2}f}{\partial x^{\beta}\partial x^{\alpha}}\frac{\partial f}{\partial x^{\alpha}}dt\nonumber \\
 &  & +4m\left(hf\right)^{\frac{2m-2-h}{h}}\frac{\partial\phi}{\partial x^{\beta}}\frac{\partial^{2}f}{\partial x^{\alpha}\partial x^{\beta}}\frac{\partial f}{\partial x^{\alpha}}dt\nonumber \\
 &  & +2m\left(hf\right)^{\frac{2m-h-2}{h}}\phi\sum_{\alpha,\beta}\left(\frac{\partial^{2}f}{\partial x^{\alpha}\partial x^{\beta}}\right)^{2}dt.\label{eq:20130813-5}
\end{eqnarray}
This is the decomposition under measure $\mathbb{P}$ for semimartingale
$H$. But in order to get more flexibility, we need to take advantage
of an important tool in stochastic analysis-the change of measure.
Let us introduce a family of probability measures on $\left(\Omega,\mathcal{F}\right)$
depending on a parameter $\lambda\in\mathbb{R}$. Define a measure
$\mathbb{Q}$ on $\left(\Omega,\mathcal{F}\right)$ by 
\begin{equation}
\left.\frac{d\mathbb{Q}}{d\mathbb{P}}\right|_{\mathcal{F}_{t}}=\exp\left(\int_{0}^{t}\sqrt{2m}\lambda\left(hf\right)^{\frac{m-2h-1}{2h}}\frac{\partial f}{\partial x^{\beta}}dW^{\beta}-m\lambda^{2}\int_{0}^{t}\left(hf\right)^{\frac{m-2h-1}{h}}\left|\nabla f\right|^{2}dt\right).\label{eq:20130318-4}
\end{equation}
Since $\left(hf\right)^{\frac{m-2h-1}{2h}}\frac{\partial f}{\partial x^{\beta}}=u^{\frac{m-3}{2}}\frac{\partial u}{\partial x^{\beta}}$,
$\beta\in\left\{ 1,\cdots,n\right\} $ are bounded by our assumption,
Novikov's condition \cite{MR1725357} is met. So the right hand side
of (\ref{eq:20130318-4}) is a true martingale under $\mathbb{P}$,
which ensures that measure $\mathbb{Q}$ is well defined. According
to Girsanov's theorem \cite{MR1121940}, under measure $\mathbb{Q}$,
the process $\tilde{W}=\left(\tilde{W}^{1},\cdots,\tilde{W}^{n}\right)$
given by 
\[
d\tilde{W}_{t}^{\beta}=dW_{t}^{\beta}-\sqrt{2m}\lambda\left(hf\right)^{\frac{m-2h-1}{2h}}\frac{\partial f}{\partial x^{\beta}}dt
\]
 is an $n-$dimensional Brownian motion. So from (\ref{eq:20130813-5}),
we easily get the decomposition for $H$ under measure $\mathbb{Q}$,
which is
\begin{eqnarray}
dH_{t} & = & 2\sqrt{2m}\left(hf\right)^{\frac{3m-2h-3}{2h}}\phi\frac{\partial^{2}f}{\partial x^{\alpha}\partial x^{\beta}}\frac{\partial f}{\partial x^{\alpha}}d\tilde{W}_{t}^{\beta}\nonumber \\
 &  & +\sqrt{2m}\left(m-h-1\right)\left(hf\right)^{\frac{3m-4h-3}{2h}}\phi\left|\nabla f\right|^{2}\frac{\partial f}{\partial x^{\beta}}d\tilde{W}_{t}^{\beta}\nonumber \\
 &  & +\sqrt{2m}\left(hf\right)^{\frac{3m-2h-3}{2h}}\frac{\partial\phi}{\partial x^{\beta}}\left|\nabla f\right|^{2}d\tilde{W}_{t}^{\beta}\nonumber \\
 &  & -m\left(2m-h+1-2\lambda\right)\left(m-h-1\right)\left(hf\right)^{\frac{2m-3h-2}{h}}\phi\left|\nabla f\right|^{4}dt\nonumber \\
 &  & +2m\left(m-h-1+\lambda\right)\left(hf\right)^{\frac{2m-2h-2}{h}}\frac{\partial\phi}{\partial x^{\alpha}}\frac{\partial f}{\partial x^{\alpha}}\left|\nabla f\right|^{2}dt\nonumber \\
 &  & +m\left(hf\right)^{\frac{2m-h-2}{h}}\Delta\phi\left|\nabla f\right|^{2}dt\nonumber \\
 &  & -2m\left(m-1\right)\left(hf\right)^{\frac{2m-2-2h}{h}}\phi\Delta f\left|\nabla f\right|^{2}dt\nonumber \\
 &  & +4m\left(\lambda-1\right)\left(hf\right)^{\frac{2m-2h-2}{h}}\phi\frac{\partial^{2}f}{\partial x^{\alpha}\partial x^{\beta}}\frac{\partial f}{\partial x^{\alpha}}\frac{\partial f}{\partial x^{\beta}}dt\nonumber \\
 &  & +4m\left(hf\right)^{\frac{2m-2-h}{h}}\frac{\partial\phi}{\partial x^{\beta}}\frac{\partial^{2}f}{\partial x^{\alpha}\partial x^{\beta}}\frac{\partial f}{\partial x^{\alpha}}dt\nonumber \\
 &  & +2m\left(hf\right)^{\frac{2m-h-2}{h}}\phi\sum_{\alpha,\beta}\left(\frac{\partial^{2}f}{\partial x^{\alpha}\partial x^{\beta}}\right)^{2}dt.\label{eq:20120921-2}
\end{eqnarray}
While from (\ref{eq:20120906-5}), it is easy to see that under measure
$\mathbb{Q}$, 
\begin{eqnarray}
 &  & dY_{t}\nonumber \\
 & = & \sqrt{2m}\left(hf\right)^{\frac{m-1}{2h}}\frac{\partial\left(f\phi\right)}{\partial x^{\alpha}}d\tilde{W}_{t}^{\alpha}\nonumber \\
 &  & +m\left(2\lambda-m+h\right)\left(hf\right)^{\frac{m-h-1}{h}}\phi\left|\nabla f\right|^{2}dt\nonumber \\
 &  & +\Delta\phi mh^{-1}\left(hf\right)^{\frac{m-1+h}{h}}dt\nonumber \\
 &  & +2m\left(1+\lambda h^{-1}\right)\left(hf\right)^{\frac{m-1}{h}}\frac{\partial f}{\partial x^{\alpha}}\frac{\partial\phi}{\partial x^{\alpha}}dt.\label{eq:20130529-1}
\end{eqnarray}
Now we have (\ref{eq:20130529-1}) and (\ref{eq:20120921-2}) at hand,
which are the fundamental decompositions we are looking for.

\subsubsection{Gradient bound for solution $u$}

In contrast to (\ref{eq:20130320-3}) which appears in the case of
global estimate for the heat equation, (\ref{eq:20120921-2}) is more
complicated, thanks to the nonlinearity of PDE $\eqref{eq:20130306-1}$
and the introduction of cut-off function $\phi$. To make use of it,
recall that we aim to estimate the first order derivatives of $f$.
Therefore, it is reasonable to get rid off the second order derivatives
of $f$ appearing in the decomposition (\ref{eq:20120921-2}). This
is done in the following lemma. For simplicity, denote by $\mathcal{A}_{x_{0}}^{R,\epsilon}$
the annulus $B\left(x_{0},R+\epsilon\right)\setminus B\left(x_{0},R\right)$.
\begin{lem}
$H$ satisfies \textup{
\begin{eqnarray}
 &  & dH_{t}\nonumber \\
 & \geq & A_{t}d\tilde{W}_{t}\nonumber \\
 &  & -\frac{1}{2}mn\left(m-1\right)^{2}\left(hf\right)^{\frac{2m-3h-2}{h}}\phi\left|\nabla f\right|^{4}dt\nonumber \\
 &  & -m\left(\left(2m-h+1-2\lambda\right)\left(m-h-1\right)+2\left(\lambda-1\right)\left(\lambda-m\right)\right)\left(hf\right)^{\frac{2m-3h-2}{h}}\phi\left|\nabla f\right|^{4}dt\nonumber \\
 &  & +2m\left(2m-h-\lambda\right)\left(hf\right)^{\frac{2m-2h-2}{h}}\frac{\partial\phi}{\partial x^{\alpha}}\frac{\partial f}{\partial x^{\alpha}}\left|\nabla f\right|^{2}dt\nonumber \\
 &  & +m\left(hf\right)^{\frac{2m-h-2}{h}}\Delta\phi\left|\nabla f\right|^{2}dt\nonumber \\
 &  & -2m\left(hf\right)^{\frac{2m-2-h}{h}}\phi^{-1}I_{\mathcal{A}_{x_{0}}^{R,\epsilon}}\left|\nabla\phi\right|^{2}\left|\nabla f\right|^{2}dt\label{eq:20120906-4}
\end{eqnarray}
}where\textup{
\begin{eqnarray}
A & = & 2\sqrt{2m}\left(hf\right)^{\frac{3m-2h-3}{2h}}\phi\frac{\partial^{2}f}{\partial x^{\alpha}\partial x^{\beta}}\frac{\partial f}{\partial x^{\alpha}}\nonumber \\
 &  & +\sqrt{2m}\left(m-h-1\right)\left(hf\right)^{\frac{3m-4h-3}{2h}}\phi\left|\nabla f\right|^{2}\frac{\partial f}{\partial x^{\beta}}\nonumber \\
 &  & +\sqrt{2m}\left(hf\right)^{\frac{3m-2h-3}{2h}}\frac{\partial\phi}{\partial x^{\beta}}\left|\nabla f\right|^{2}.\label{eq:2012-12-4-2}
\end{eqnarray}
}\end{lem}
\begin{proof}
From (\ref{eq:20120921-2}), let us write 
\[
dH_{t}\triangleq A_{t}d\tilde{W}_{t}+\left(B+C+D\right)dt
\]
where 
\begin{eqnarray*}
A & = & 2\sqrt{2m}\left(hf\right)^{\frac{3m-2h-3}{2h}}\phi\frac{\partial^{2}f}{\partial x^{\alpha}\partial x^{\beta}}\frac{\partial f}{\partial x^{\alpha}}\\
 &  & +\sqrt{2m}\left(m-h-1\right)\left(hf\right)^{\frac{3m-4h-3}{2h}}\phi\left|\nabla f\right|^{2}\frac{\partial f}{\partial x^{\beta}}\\
 &  & +\sqrt{2m}\left(hf\right)^{\frac{3m-2h-3}{2h}}\frac{\partial\phi}{\partial x^{\beta}}\left|\nabla f\right|^{2}
\end{eqnarray*}
and by Cauchy-Schwartz inequality, 
\begin{eqnarray*}
B & = & \sum_{\alpha}\left(2m\left(hf\right)^{\frac{2m-h-2}{h}}\phi\left(\frac{\partial^{2}f}{\partial x^{\alpha}\partial x^{\alpha}}\right)^{2}\right.\\
 &  & -2m\left(m-1\right)\left(hf\right)^{\frac{2m-2-2h}{h}}\phi\frac{\partial^{2}f}{\partial x^{\alpha}\partial x^{\alpha}}\left|\nabla f\right|^{2}\\
 &  & +4m\left(\lambda-1\right)\left(hf\right)^{\frac{2m-2h-2}{h}}\phi\frac{\partial^{2}f}{\partial x^{\alpha}\partial x^{\alpha}}\frac{\partial f}{\partial x^{\alpha}}\frac{\partial f}{\partial x^{\alpha}}\\
 &  & \left.+4m\left(hf\right)^{\frac{2m-2-h}{h}}\frac{\partial\phi}{\partial x^{\alpha}}\frac{\partial^{2}f}{\partial x^{\alpha}\partial x^{\alpha}}\frac{\partial f}{\partial x^{\alpha}}\right)\\
 & \geq & -2m\left(hf\right)^{\frac{2m-2-h}{h}}\phi^{-1}I_{\mathcal{A}_{x_{0}}^{R,\epsilon}}\sum_{\alpha}\left(\frac{\partial\phi}{\partial x^{\alpha}}\right)^{2}\left(\frac{\partial f}{\partial x^{\alpha}}\right)^{2}\\
 &  & -2m\left(\lambda-1\right)^{2}\left(hf\right)^{\frac{2m-3h-2}{h}}\phi\sum_{\alpha}\left(\frac{\partial f}{\partial x^{\alpha}}\right)^{4}\\
 &  & -4m\left(\lambda-1\right)\left(hf\right)^{\frac{2m-2h-2}{h}}\sum_{\alpha}\frac{\partial\phi}{\partial x^{\alpha}}\left(\frac{\partial f}{\partial x^{\alpha}}\right)^{3}\\
 &  & +2m\left(m-1\right)\left(hf\right)^{\frac{2m-2-2h}{h}}\left|\nabla f\right|^{2}\frac{\partial\phi}{\partial x^{\alpha}}\frac{\partial f}{\partial x^{\alpha}}\\
 &  & +\left(\lambda-1-\frac{1}{4}n\left(m-1\right)\right)2m\left(m-1\right)\left(hf\right)^{\frac{2m-2-3h}{h}}\phi\left|\nabla f\right|^{4},
\end{eqnarray*}
 
\begin{eqnarray*}
C & = & \sum_{\alpha\neq\beta}\left(2m\left(hf\right)^{\frac{2m-h-2}{h}}\phi\left(\frac{\partial^{2}f}{\partial x^{\alpha}\partial x^{\beta}}\right)^{2}\right.\\
 &  & +4m\left(\lambda-1\right)\left(hf\right)^{\frac{2m-2h-2}{h}}\phi\frac{\partial^{2}f}{\partial x^{\alpha}\partial x^{\beta}}\frac{\partial f}{\partial x^{\alpha}}\frac{\partial f}{\partial x^{\beta}}\\
 &  & \left.+4m\left(hf\right)^{\frac{2m-2-h}{h}}\frac{\partial\phi}{\partial x^{\beta}}\frac{\partial^{2}f}{\partial x^{\alpha}\partial x^{\beta}}\frac{\partial f}{\partial x^{\alpha}}\right)\\
 & \geq & \sum_{\alpha\neq\beta}\left(-2m\left(\lambda-1\right)^{2}\left(hf\right)^{\frac{2m-3h-2}{h}}\phi\left(\frac{\partial f}{\partial x^{\alpha}}\frac{\partial f}{\partial x^{\beta}}\right)^{2}\right.\\
 &  & -2m\left(hf\right)^{\frac{2m-2-h}{h}}\phi^{-1}I_{\mathcal{A}_{x_{0}}^{R,\epsilon}}\left(\frac{\partial\phi}{\partial x^{\beta}}\frac{\partial f}{\partial x^{\alpha}}\right)^{2}\\
 &  & \left.-4m\left(\lambda-1\right)\left(hf\right)^{\frac{2m-2h-2}{h}}\frac{\partial\phi}{\partial x^{\beta}}\frac{\partial f}{\partial x^{\alpha}}\frac{\partial f}{\partial x^{\alpha}}\frac{\partial f}{\partial x^{\beta}}\right)
\end{eqnarray*}
and
\begin{eqnarray*}
D & = & -m\left(2m-h+1-2\lambda\right)\left(m-h-1\right)\left(hf\right)^{\frac{2m-3h-2}{h}}\phi\left|\nabla f\right|^{4}\\
 &  & +m\left(hf\right)^{\frac{2m-h-2}{h}}\Delta\phi\left|\nabla f\right|^{2}\\
 &  & +2m\left(m-h-1+\lambda\right)\left(hf\right)^{\frac{2m-2h-2}{h}}\frac{\partial\phi}{\partial x^{\alpha}}\frac{\partial f}{\partial x^{\alpha}}\left|\nabla f\right|^{2}.
\end{eqnarray*}
Then adding these three inequalities gives us (\ref{eq:20120906-4}). 
\end{proof}
To deal with the finite variation part in (\ref{eq:20120906-4}),
let us assume that 
\begin{equation}
\frac{1}{2}n\left(m-1\right)^{2}+\left(2m-h+1-2\lambda\right)\left(m-h-1\right)+2\left(\lambda-1\right)\left(\lambda-m\right)<0.\label{eq:20120906-3}
\end{equation}
This allows us to get the following estimate. 
\begin{lem}
Under assumption \textup{(\ref{eq:20120906-3})}, we have
\begin{equation}
E^{\mathbb{Q}}\left[\int_{0}^{T}H_{s}ds\right]\geq H_{0}T-\left(\frac{1}{2}\left|2m-h-\lambda\right|l^{-3}+\frac{3}{2}l^{-1}\right)mM_{1}L^{2}\frac{T^{2}}{2}\label{eq:20120911-1}
\end{equation}
where\textup{
\begin{equation}
l=\frac{\frac{1}{2}n\left(m-1\right)^{2}+\left(2m-h+1-2\lambda\right)\left(m-h-1\right)+2\left(\lambda-1\right)\left(\lambda-m\right)}{-\frac{3}{2}\left|2m-h-\lambda\right|-\frac{3}{2}}\label{eq:20130321-5}
\end{equation}
}and\textup{ 
\[
M_{1}=\sup_{\left[0,T\right]\times\mathcal{A}_{x_{0}}^{R,\epsilon}}u^{2m+h-2}\left(t,y\right).
\]
}\end{lem}
\begin{proof}
By Cauchy-Schwartz inequality, for any positive $l_{1}$, $l_{2}$,
$l_{3}$ and $l_{4}$, 
\begin{eqnarray*}
 &  & \left|\left(hf\right)^{\frac{2m-2h-2}{h}}\frac{\partial\phi}{\partial x^{\alpha}}\frac{\partial f}{\partial x^{\alpha}}\left|\nabla f\right|^{2}\right|\\
 & \leq & \frac{1}{2}l_{1}\left(hf\right)^{\frac{2m-3h-2}{h}}\phi\left|\nabla f\right|^{4}+\frac{1}{2}l_{1}^{-1}\left(hf\right)^{\frac{2m-h-2}{h}}\phi^{-1}\left|\nabla\phi\right|^{2}\left|\nabla f\right|^{2}\\
 & \leq & \frac{1}{2}l_{1}\left(hf\right)^{\frac{2m-3h-2}{h}}\phi\left|\nabla f\right|^{4}+\frac{1}{4}l_{1}^{-1}l_{2}\left(hf\right)^{\frac{2m-3h-2}{h}}\phi\left|\nabla f\right|^{4}\\
 &  & +\frac{1}{4}l_{1}^{-1}l_{2}^{-1}\left(hf\right)^{\frac{2m+h-2}{h}}\phi^{-3}\left|\nabla\phi\right|^{4},
\end{eqnarray*}
 
\begin{eqnarray*}
 &  & \left|\left(hf\right)^{\frac{2m-h-2}{h}}\phi^{-1}\left|\nabla\phi\right|^{2}\left|\nabla f\right|^{2}\right|\\
 & \leq & \frac{1}{2}l_{3}\left(hf\right)^{\frac{2m-3h-2}{h}}\phi\left|\nabla f\right|^{4}+\frac{1}{2}l_{3}^{-1}\left(hf\right)^{\frac{2m+h-2}{h}}\phi^{-3}\left|\nabla\phi\right|^{4}
\end{eqnarray*}
and 
\begin{eqnarray*}
 &  & \left|\left(hf\right)^{\frac{2m-h-2}{h}}\Delta\phi\left|\nabla f\right|^{2}\right|\\
 & \leq & \frac{1}{2}l_{4}\left(hf\right)^{\frac{2m-3h-2}{h}}\phi\left|\nabla f\right|^{4}+\frac{1}{2}l_{4}^{-1}\left(hf\right)^{\frac{2m+h-2}{h}}\phi^{-1}\left|\Delta\phi\right|^{2}.
\end{eqnarray*}
Plugging them into (\ref{eq:20120906-4}) in the above lemma yields
\begin{eqnarray*}
 &  & dH_{t}\\
 & \geq & A_{t}d\tilde{W}_{t}\\
 &  & -\frac{1}{2}mn\left(m-1\right)^{2}\left(hf\right)^{\frac{2m-3h-2}{h}}\phi\left|\nabla f\right|^{4}dt\\
 &  & -m\left(\left(2m-h+1-2\lambda\right)\left(m-h-1\right)+2\left(\lambda-1\right)\left(\lambda-m\right)\right)\left(hf\right)^{\frac{2m-3h-2}{h}}\phi\left|\nabla f\right|^{4}dt\\
 &  & -2m\left|2m-h-\lambda\right|\left(\frac{1}{2}l_{1}+\frac{1}{4}l_{1}^{-1}l_{2}\right)\left(hf\right)^{\frac{2m-3h-2}{h}}\phi\left|\nabla f\right|^{4}dt\\
 &  & -ml_{3}\left(hf\right)^{\frac{2m-3h-2}{h}}\phi\left|\nabla f\right|^{4}dt\\
 &  & -\frac{1}{2}ml_{4}\left(hf\right)^{\frac{2m-3h-2}{h}}\phi\left|\nabla f\right|^{4}dt\\
 &  & -\frac{1}{2}ml_{4}^{-1}\left(hf\right)^{\frac{2m+h-2}{h}}\phi^{-1}\left|\Delta\phi\right|^{2}dt\\
 &  & -\frac{1}{2}m\left|2m-h-\lambda\right|l_{1}^{-1}l_{2}^{-1}\left(hf\right)^{\frac{2m+h-2}{h}}\phi^{-3}\left|\nabla\phi\right|^{4}dt\\
 &  & -ml_{3}^{-1}\left(hf\right)^{\frac{2m+h-2}{h}}\phi^{-3}\left|\nabla\phi\right|^{4}dt.
\end{eqnarray*}
For simplicity, set $l_{1}^{2}=l_{2}=l_{3}^{2}=l_{4}^{2}=l^{2}$.
Then we have
\begin{eqnarray*}
 &  & dH_{t}\\
 & \geq & A_{t}d\tilde{W}_{t}\\
 &  & -\frac{1}{2}mn\left(m-1\right)^{2}\left(hf\right)^{\frac{2m-3h-2}{h}}\phi\left|\nabla f\right|^{4}dt\\
 &  & -m\left(\left(2m-h+1-2\lambda\right)\left(m-h-1\right)+2\left(\lambda-1\right)\left(\lambda-m\right)\right)\left(hf\right)^{\frac{2m-3h-2}{h}}\phi\left|\nabla f\right|^{4}dt\\
 &  & -\frac{3}{2}m\left|2m-h-\lambda\right|l\left(hf\right)^{\frac{2m-3h-2}{h}}\phi\left|\nabla f\right|^{4}dt\\
 &  & -\frac{3}{2}ml\left(hf\right)^{\frac{2m-3h-2}{h}}\phi\left|\nabla f\right|^{4}dt\\
 &  & -\frac{1}{2}m\left(\left|2m-h-\lambda\right|l^{-3}+2l^{-1}\right)\left(hf\right)^{\frac{2m+h-2}{h}}\phi^{-3}\left|\nabla\phi\right|^{4}dt\\
 &  & -\frac{1}{2}ml^{-1}\left(hf\right)^{\frac{2m+h-2}{h}}\phi^{-1}\left|\Delta\phi\right|^{2}dt.
\end{eqnarray*}
By estimates for the cut-off function $\phi$ in  (\ref{eq:2012-11-28-2})
and (\ref{eq:2012-11-28-3}),
\begin{eqnarray*}
 &  & dH_{t}\\
 & \geq & A_{t}d\tilde{W}_{t}\\
 &  & -\frac{1}{2}mn\left(m-1\right)^{2}\left(hf\right)^{\frac{2m-3h-2}{h}}\phi\left|\nabla f\right|^{4}dt\\
 &  & -m\left(\left(2m-h+1-2\lambda\right)\left(m-h-1\right)+2\left(\lambda-1\right)\left(\lambda-m\right)\right)\left(hf\right)^{\frac{2m-3h-2}{h}}\phi\left|\nabla f\right|^{4}dt\\
 &  & -\frac{3}{2}m\left|2m-h-\lambda\right|l\left(hf\right)^{\frac{2m-3h-2}{h}}\phi\left|\nabla f\right|^{4}dt\\
 &  & -\frac{3}{2}ml\left(hf\right)^{\frac{2m-3h-2}{h}}\phi\left|\nabla f\right|^{4}dt\\
 &  & -\frac{1}{2}m\left(\left|2m-h-\lambda\right|l^{-3}+3l^{-1}\right)\left(hf\right)^{\frac{2m+h-2}{h}}I_{\mathcal{A}_{x_{0}}^{R,\epsilon}}L^{2}dt.
\end{eqnarray*}
 Now, by assumption (\ref{eq:20120906-3}), we can choose a positive
$l$ small enough such that
\begin{eqnarray*}
 &  & -\frac{1}{2}mn\left(m-1\right)^{2}-m\left(\left(2m-h+1-2\lambda\right)\left(m-h-1\right)+2\left(\lambda-1\right)\left(\lambda-m\right)\right)\\
 & = & \left(\frac{3}{2}\left|2m-h-\lambda\right|+\frac{3}{2}\right)ml.
\end{eqnarray*}
Hence we have for any $t\in\left[0,T\right]$, 
\begin{eqnarray*}
dH_{t} & \geq & A_{t}d\tilde{W}_{t}-\left(\frac{1}{2}\left|2m-h-\lambda\right|l^{-3}+\frac{3}{2}l^{-1}\right)mM_{1}L^{2}dt
\end{eqnarray*}
 where
\[
M_{1}=\sup_{\begin{array}{c}
\left[0,T\right]\times\mathcal{A}_{x_{0}}^{R,\epsilon}\end{array}}\left|\left(hf\right)^{\frac{2m+h-2}{h}}\left(t,y\right)\right|=\sup_{\left[0,T\right]\times\mathcal{A}_{x_{0}}^{R,\epsilon}}u^{2m+h-2}\left(t,y\right).
\]
 This gives
\begin{equation}
H_{s}-H_{0}\geq\int_{0}^{s}A_{t}d\tilde{W}_{t}-\left(\frac{1}{2}\left|2m-h-\lambda\right|l^{-3}+\frac{3}{2}l^{-1}\right)mM_{1}L^{2}s\label{eq:2012-12-4-3}
\end{equation}
for any $s\in\left[0,T\right]$. According to our assumption on $u$,
$A$, defined by (\ref{eq:2012-12-4-2}), is bounded. So $\int_{0}^{\cdot}A_{s}d\tilde{W}_{s}$
is a true martingale under $\mathbb{Q}$. Therefore we can take expectation
on both sides of (\ref{eq:2012-12-4-3}) and obtain
\begin{equation}
E^{\mathbb{Q}}\left[H_{s}\right]\geq H_{0}-\left(\frac{1}{2}\left|2m-h-\lambda\right|l^{-3}+\frac{3}{2}l^{-1}\right)mM_{1}L^{2}s.\label{eq:20131004-3}
\end{equation}
 Integrating both sides with respect to $s$ on $\left[0,T\right]$
yields (\ref{eq:20120911-1}). 
\end{proof}
On the other hand, by using (\ref{eq:20120906-5}), we have an upper
bound on the left hand side of (\ref{eq:20120911-1}).
\begin{lem}
Assume\textup{ 
\begin{equation}
2\lambda-m+h>0.\label{eq:20120911-2-1}
\end{equation}
}Then for any\textup{ $\rho>0$ }such that\textup{ 
\begin{equation}
2\lambda-m+h-\rho\left|1+\lambda h^{-1}\right|>0,\label{eq:20130321-4}
\end{equation}
}we have
\begin{equation}
E^{\mathbb{Q}}\left[\int_{0}^{T}H_{t}dt\right]\leq\frac{E^{\mathbb{Q}}\left[Y_{T}\right]-Y_{0}+\left(\left|h^{-1}\right|+\rho^{-1}\left|1+\lambda h^{-1}\right|\right)mLM_{2}T}{m\left(2\lambda-m+h-\rho\left|1+\lambda h^{-1}\right|\right)}\label{eq:20120921-3}
\end{equation}
with\textup{
\[
M_{2}=\sup_{\left[0,T\right]\times\mathcal{A}_{x_{0}}^{R,\epsilon}}u^{m+h-1}\left(t,y\right).
\]
}\end{lem}
\begin{proof}
Recall (\ref{eq:20130529-1}), which is the fundamental decomposition
of $Y$ under measure $\mathbb{Q}$. Then by property (\ref{eq:2012-11-28-2})
of cut-off function $\phi$ and Cauchy-Schwartz inequality, for any
$\rho>0$, 
\begin{eqnarray}
 &  & dY_{t}\nonumber \\
 & \geq & \sqrt{2m}\left(hf\right)^{\frac{m-1}{2h}}\frac{\partial\left(f\phi\right)}{\partial x^{\alpha}}d\tilde{W}_{t}^{\alpha}\nonumber \\
 &  & +m\left(2\lambda-m+h\right)\left(hf\right)^{\frac{m-h-1}{h}}\phi\left|\nabla f\right|^{2}dt\nonumber \\
 &  & -L\phi^{\frac{1}{2}}I_{\mathcal{A}_{x_{0}}^{R,\epsilon}}m\left|h^{-1}\right|\left(hf\right)^{\frac{m-1+h}{h}}dt\nonumber \\
 &  & -\rho m\left|1+\lambda h^{-1}\right|\left(hf\right)^{\frac{m-1-h}{h}}\phi\left|\nabla f\right|^{2}dt\nonumber \\
 &  & -\rho^{-1}m\left|1+\lambda h^{-1}\right|\left(hf\right)^{\frac{m-1+h}{h}}I_{\mathcal{A}_{x_{0}}^{R,\epsilon}}\phi^{-1}\left|\nabla\phi\right|^{2}dt\nonumber \\
 & \geq & \sqrt{2m}\left(hf\right)^{\frac{m-1}{2h}}\frac{\partial\left(f\phi\right)}{\partial x^{\alpha}}d\tilde{W}_{t}^{\alpha}\nonumber \\
 &  & +m\left(2\lambda-m+h-\rho\left|1+\lambda h^{-1}\right|\right)H_{t}dt\nonumber \\
 &  & -\left(\left|h^{-1}\right|+\rho^{-1}\left|1+\lambda h^{-1}\right|\right)mLM_{2}dt\label{eq:20130501-1}
\end{eqnarray}
where
\[
M_{2}=\sup_{\left[0,T\right]\times\mathcal{A}_{x_{0}}^{R,\epsilon}}\left(hf\right)^{\frac{m+h-1}{h}}\left(t,y\right)=\sup_{\left[0,T\right]\times\mathcal{A}_{x_{0}}^{R,\epsilon}}u^{m+h-1}\left(t,y\right).
\]
According to our assumption on $u$, $\left(hf\right)^{\frac{m-1}{2h}}\frac{\partial\left(f\phi\right)}{\partial x^{\alpha}}$,
$\alpha\in\left\{ 1,\cdots,n\right\} $ are all bounded, which ensures
that $\int_{0}^{\cdot}\left(hf\right)^{\frac{m-1}{2h}}\frac{\partial\left(f\phi\right)}{\partial x^{\alpha}}d\tilde{W}_{t}^{\alpha}$
is a $\mathbb{Q}$ martingale. Hence
\begin{eqnarray*}
 &  & E^{\mathbb{Q}}\left[Y_{T-t_{1}}\right]-Y_{0}\\
 & \geq & m\left(2\lambda-m+h-\rho\left|1+\lambda h^{-1}\right|\right)E^{\mathbb{Q}}\left[\int_{0}^{T-t_{1}}H_{t}dt\right]\\
 &  & -\left(\left|h^{-1}\right|+\rho^{-1}\left|1+\lambda h^{-1}\right|\right)mLM_{2}\left(T-t_{1}\right).
\end{eqnarray*}
Dividing both sides by $m\left(2\lambda-m+h-\rho\left|1+\lambda h^{-1}\right|\right)$
yields (\ref{eq:20120921-3}).
\end{proof}
Now, based on (\ref{eq:20120911-1}) and (\ref{eq:20120921-3}), we
are ready to obtain a gradient bound for $u$. 
\begin{prop}
Let the assumptions in Lemma 8 and Lemma 9 be satisfied, then for
any\textup{ $\left(T,x\right)\in\left(0,t_{1}\right)\times B\left(x_{0},R\right)$},\textup{
\begin{eqnarray}
 &  & u^{m+h-3}\left|\nabla u\right|^{2}\left(T,x\right)\nonumber \\
 & \leq & \frac{E^{\mathbb{Q}}\left[\phi u_{0}^{h}\left(X_{T}\right)\right]-u^{h}\left(T,x\right)}{hm\left(2\lambda-m+h-\rho\left|1+\lambda h^{-1}\right|\right)T}+\frac{\left(\left|h^{-1}\right|+\rho^{-1}\left|1+\lambda h^{-1}\right|\right)LM_{2}}{2\lambda-m+h-\rho\left|1+\lambda h^{-1}\right|}\nonumber \\
 &  & +\left(\frac{1}{2}\left|2m-h-\lambda\right|l^{-3}+\frac{3}{2}l^{-1}\right)mM_{1}L^{2}\frac{T}{2}.\label{eq:20120911-4}
\end{eqnarray}
}\end{prop}
\begin{proof}
Combining (\ref{eq:20120911-1}) with (\ref{eq:20120921-3}), we have
\begin{eqnarray*}
 &  & H_{0}\\
 & \leq & \frac{E^{\mathbb{Q}}\left[Y_{T}\right]-Y_{0}+\left(\left|h^{-1}\right|+\rho^{-1}\left|1+\lambda h^{-1}\right|\right)mLM_{2}T}{m\left(2\lambda-m+h-\rho\left|1+\lambda h^{-1}\right|\right)T}\\
 &  & +\left(\frac{1}{2}\left|2m-h-\lambda\right|l^{-3}+\frac{3}{2}l^{-1}\right)mM_{1}L^{2}\frac{T}{2}.
\end{eqnarray*}
Then we arrive at (\ref{eq:20120911-4}) by noting that 
\[
H_{0}=\left(hf\right)^{\frac{m-h-1}{h}}\phi\left|\nabla f\right|^{2}\left(T,x\right),
\]
\[
Y_{0}=\phi f\left(T,x\right)=\phi\frac{u^{h}}{h}\left(T,x\right)
\]
 and 
\[
E^{\mathbb{Q}}\left[Y_{T}\right]=E^{\mathbb{Q}}\left[\phi\frac{u_{0}^{h}}{h}\left(X_{T}\right)\right].
\]

\end{proof}
Let us revise the assumptions we have made, that is, (\ref{eq:20120906-3})
and (\ref{eq:20120911-2-1}). They impose a restriction on the choice
of $m$. Let us look at (\ref{eq:20120906-3}) first. It is equivalent
to
\begin{equation}
2\lambda^{2}-2\left(2m-h\right)\lambda+\frac{1}{2}n\left(m-1\right)^{2}+\left(2m-h+1\right)\left(m-h-1\right)+2m<0\label{eq:20130321-1}
\end{equation}
To let the left hand side attain the minimum, we take $\lambda=\frac{2m-h}{2}$.
Then (\ref{eq:20130321-1}) is reduced to
\begin{equation}
h^{2}-2mh+n\left(m-1\right)^{2}+2\left(m-1\right)<0.\label{eq:2012-12-4-1}
\end{equation}
Notice that after setting $\lambda=\frac{2m-h}{2}$, condition (\ref{eq:20120911-2-1})
is automatically satisfied. Moreover, (\ref{eq:2012-12-4-1}) holds
if and only if $m\in\left(1-\frac{1}{\sqrt{n-1}},1+\frac{1}{\sqrt{n-1}}\right)$,
and 
\[
h\in\left(m-\sqrt{1-\left(n-1\right)\left(m-1\right)^{2}},m+\sqrt{1-\left(n-1\right)\left(m-1\right)^{2}}\right).
\]
 Then we arrive at Theorem 1 by substituting $\lambda=\frac{2m-h}{2}$
into (\ref{eq:20130321-4}) and (\ref{eq:20130321-5}), and noting
that $E^{\mathbb{Q}}\left[\phi u_{0}^{h}\left(X_{T}\right)\right]\leq\left\Vert I_{B\left(x_{0},R+\epsilon\right)}u_{0}^{h}\right\Vert _{\infty}.$

\subsubsection{Gradient bound for pressure variable $v=\frac{u^{m-1}}{m-1}$}

Since most of the existing results are in terms of $\frac{\left|\nabla v\right|^{2}}{v}$,
we are tempted to work out a bound on it with our method. Notice that
$\frac{\left|\nabla v\right|^{2}}{v}=\left(m-1\right)u^{m-3}\left|\nabla u\right|^{2}$,
which corresponds to the LHS of (\ref{eq:20130528-1}) in Theorem
1 when $h=0$. However, the RHS of (\ref{eq:20130528-1}) will explode
when $h\rightarrow0$. Therefore, we can not use the result in Theorem
1 directly. 

Let us modify the estimate we obtained from the decomposition of $Y$,
so that the bound does not explode when $h\rightarrow0$. The idea
is to decompose $Y_{t}-\frac{\phi}{h}\left(X_{t}\right)$ instead
of $Y_{t}$. This is because 
\[
Y_{t}-\frac{\phi}{h}\left(X_{t}\right)=\phi\frac{u^{h}-1}{h}\rightarrow\phi\log u,
\]
which is finite.
\begin{lem}
Assume\textup{ 
\begin{equation}
2\lambda-m+h>0.\label{eq:20120911-2-1-1}
\end{equation}
}Then for any\textup{ $\rho>0$ }such that\textup{ 
\begin{equation}
2\lambda-m+h-\rho>0,\label{eq:20130321-4-1}
\end{equation}
}we have
\begin{equation}
E^{\mathbb{Q}}\left[\int_{0}^{T}H_{t}dt\right]\leq\frac{E^{\mathbb{Q}}\left[Y_{T}-\frac{\phi}{h}\left(X_{T}\right)\right]-Y_{0}+\frac{\phi}{h}\left(x\right)+Lm\rho^{-1}M_{3}T}{m\left(2\lambda-m+h-\rho\right)}\label{eq:20130606-3}
\end{equation}
\textup{ }where
\[
M_{3}=\sup_{\left[0,T\right]\times\mathcal{A}_{x_{0}}^{R,\epsilon}}\left(\rho\left|\frac{u^{h}-1}{h}\right|u^{h}+\left(\lambda\left(\frac{u^{h}-1}{h}\right)+u^{h}\right)^{2}\right)u^{m-h-1}\left(t,y\right).
\]
\end{lem}
\begin{proof}
From (\ref{eq:20130529-1}), 
\begin{eqnarray}
 &  & d\left(Y_{t}-\frac{\phi}{h}\left(X_{t}\right)\right)\nonumber \\
 & = & \sqrt{2m}\left(hf\right)^{\frac{m-1}{2h}}\left(\phi\frac{\partial f}{\partial x^{\alpha}}+\left(f-\frac{1}{h}\right)\frac{\partial\phi}{\partial x^{\alpha}}\right)d\tilde{W}_{t}^{\alpha}\nonumber \\
 &  & +m\left(2\lambda-m+h\right)\left(hf\right)^{\frac{m-h-1}{h}}\phi\left|\nabla f\right|^{2}dt\nonumber \\
 &  & +\Delta\phi m\left(f-\frac{1}{h}\right)\left(hf\right)^{\frac{m-1}{h}}dt\nonumber \\
 &  & +2m\left(\lambda\left(f-\frac{1}{h}\right)+hf\right)\left(hf\right)^{\frac{m-h-1}{h}}\frac{\partial f}{\partial x^{\alpha}}\frac{\partial\phi}{\partial x^{\alpha}}dt.\label{eq:20130607-1}
\end{eqnarray}
Then by Cauchy-Schwartz inequality, for any $\rho>0$,
\begin{eqnarray*}
 &  & d\left(Y_{t}-\frac{\phi}{h}\left(X_{t}\right)\right)\\
 & \geq & \sqrt{2m}\left(hf\right)^{\frac{m-1}{2h}}\left(\phi\frac{\partial f}{\partial x^{\alpha}}+\left(f-\frac{1}{h}\right)\frac{\partial\phi}{\partial x^{\alpha}}\right)d\tilde{W}_{t}^{\alpha}\\
 &  & +m\left(2\lambda-m+h\right)\left(hf\right)^{\frac{m-h-1}{h}}\phi\left|\nabla f\right|^{2}dt\\
 &  & -L\phi^{\frac{1}{2}}I_{\mathcal{A}_{x_{0}}^{R,\epsilon}}m\left|f-\frac{1}{h}\right|\left(hf\right)^{\frac{m-1}{h}}dt\\
 &  & -\rho m\left(hf\right)^{\frac{m-h-1}{h}}\phi\left|\nabla f\right|^{2}dt\\
 &  & -\rho^{-1}m\left(\lambda\left(f-\frac{1}{h}\right)+hf\right)^{2}\left(hf\right)^{\frac{m-h-1}{h}}I_{\mathcal{A}_{x_{0}}^{R,\epsilon}}\phi^{-1}\left|\nabla\phi\right|^{2}dt\\
 & \geq & \sqrt{2m}\left(hf\right)^{\frac{m-1}{2h}}\left(\phi\frac{\partial f}{\partial x^{\alpha}}+\left(f-\frac{1}{h}\right)\frac{\partial\phi}{\partial x^{\alpha}}\right)d\tilde{W}_{t}^{\alpha}\\
 &  & +m\left(2\lambda-m+h-\rho\right)H_{t}dt\\
 &  & -Lm\rho^{-1}\left(\rho\left|f-\frac{1}{h}\right|hf+\left(\lambda\left(f-\frac{1}{h}\right)+hf\right)^{2}\right)\left(hf\right)^{\frac{m-h-1}{h}}dt.
\end{eqnarray*}
Let
\[
M_{3}=\sup_{\left[0,T\right]\times\mathcal{A}_{x_{0}}^{R,\epsilon}}\left(\rho\left|\frac{u^{h}-1}{h}\right|u^{h}+\left(\lambda\left(\frac{u^{h}-1}{h}\right)+u^{h}\right)^{2}\right)u^{m-h-1}\left(t,y\right).
\]
According to our assumption on $u$, $\left(hf\right)^{\frac{m-1}{2h}}\left(\phi\frac{\partial f}{\partial x^{\alpha}}+\left(f-\frac{1}{h}\right)\frac{\partial\phi}{\partial x^{\alpha}}\right)$,
$\alpha\in\left\{ 1,\cdots,n\right\} $ are all bounded, which ensures
that the local martingale part of $Y_{t}-\frac{\phi}{h}\left(X_{t}\right)$
is a true $\mathbb{Q}$ martingale. Hence
\begin{eqnarray*}
 &  & E^{\mathbb{Q}}\left[Y_{T}-\frac{\phi}{h}\left(X_{T}\right)\right]-Y_{0}+\frac{\phi}{h}\left(x\right)\\
 & \geq & m\left(2\lambda-m+h-\rho\right)E^{\mathbb{Q}}\left[\int_{0}^{T}H_{t}dt\right]-Lm\rho^{-1}M_{3}T,
\end{eqnarray*}
which completes the proof.
\end{proof}
By combining (\ref{eq:20130606-3}) with (\ref{eq:20120911-1}) in
the way as in the proof of Theorem 1 and also taking $\lambda=\frac{2m-h}{2}$,
we are able to get the following gradient bound.
\begin{thm}
If $u$ is a positive and bounded solution to \textup{(\ref{eq:20130306-1})}
on \textup{$\left(0,t_{1}\right)\times B\left(x_{0},R+\epsilon\right)$}
with $m\in\left(1-\frac{1}{\sqrt{n-1}},1+\frac{1}{\sqrt{n-1}}\right)$,
then for any 
\[
h\in\left(m-\sqrt{1-\left(n-1\right)\left(m-1\right)^{2}},m+\sqrt{1-\left(n-1\right)\left(m-1\right)^{2}}\right),
\]
 $\rho\in\left(0,m\right)$, and $\left(T,x\right)\in\left(0,t_{1}\right)\times B\left(x_{0},R\right)$,
we have 
\begin{eqnarray*}
 &  & u^{m+h-3}\left|\nabla u\right|^{2}\left(T,x\right)\\
 & \leq & \frac{E^{\mathbb{Q}}\left[\phi\left(u_{0}^{h}-1\right)\left(X_{T}\right)\right]-\left(u^{h}\left(T,x\right)-1\right)}{hm\left(m-\rho\right)T}\\
 &  & +\frac{LM_{3}}{\rho\left(m-\rho\right)}+\left(\frac{1}{4}\left|2m-h\right|l^{-3}+\frac{3}{2}l^{-1}\right)mM_{1}L^{2}\frac{T}{2}
\end{eqnarray*}
where 
\[
l=\frac{-\frac{1}{2}n\left(m-1\right)^{2}+\left(h-1\right)\left(m-1\right)-\frac{1}{2}h^{2}+h}{\frac{3}{4}\left|2m-h\right|+\frac{3}{2}},\ M_{1}=\sup_{\left[0,T\right]\times\mathcal{A}_{x_{0}}^{R,\epsilon}}u^{2m+h-2},
\]
\textup{
\[
M_{3}=\sup_{\left[0,T\right]\times\mathcal{A}_{x_{0}}^{R,\epsilon}}\left(\rho\left|\frac{u^{h}-1}{h}\right|u^{h}+\left(\frac{2m-h}{2}\left(\frac{u^{h}-1}{h}\right)+u^{h}\right)^{2}\right)u^{m-h-1}\left(t,y\right),
\]
}and $L$ is a constant depending on the cut-off function $\phi$
by (\ref{eq:2012-11-28-2}) and (\ref{eq:2012-11-28-3}). 
\end{thm}
By letting $h\rightarrow0$, we immediately get Corollary 2.

\subsection{Global gradient estimate}

In this section, we consider a positive and bounded solution $u$
to (\ref{eq:20130306-1}) on $\left(0,\infty\right)\times\mathbb{R}^{n}$.
First of all, as a direct consequence of the main result in the last
section, we can get a global gradient bound from the local bound by
letting the radius of the local ball tend to infinity.

\subsubsection{From local bound to global bound}

Since $u$ is positive and bounded on $\left(0,\infty\right)\times\mathbb{R}^{n}$,
the local result in Theorem 1 holds for any $t_{1}$, $R$ and $\epsilon$.
By substituting $\epsilon=R$ into (\ref{eq:20130528-1}), and then
taking $R\rightarrow\infty$ on both sides, we obtain
\[
u^{m+h-3}\left|\nabla u\right|^{2}\left(T,x\right)\leq\frac{E^{\mathbb{Q}}\left[u_{0}^{h}\left(X_{T}\right)\right]-u^{h}\left(T,x\right)}{Thm\left(m-\rho\left|\frac{2m+h}{2h}\right|\right)}
\]
where we have used the fact that when $\epsilon=R$ 
\[
\lim_{R\rightarrow\infty}L=\lim_{R\rightarrow\infty}\frac{8}{R^{2}}\left(\frac{nR}{2R+R}+\frac{8\left(2R\right)^{2}}{\left(3R\right)^{2}}\right)=0\mbox{ and }\lim_{R\rightarrow\infty}\phi=1.
\]
By letting $\rho=0$, we have 
\begin{equation}
u^{m+h-3}\left|\nabla u\right|^{2}\left(T,x\right)\leq\frac{E^{\mathbb{Q}}\left[u_{0}^{h}\left(X_{T}\right)\right]-u^{h}\left(T,x\right)}{Thm^{2}}.\label{eq:20130605-1}
\end{equation}
 Note that in the local case, due to the existence of $l^{-1}$ in
(\ref{eq:20130528-1}), $h$ is not allowed to touch the boundary
of the open interval 
\[
\left(m-\sqrt{1-\left(n-1\right)\left(m-1\right)^{2}},m+\sqrt{1-\left(n-1\right)\left(m-1\right)^{2}}\right).
\]
But now we can take limit of $h$ to 
\[
m-\sqrt{1-\left(n-1\right)\left(m-1\right)^{2}}\mbox{ or }m+\sqrt{1-\left(n-1\right)\left(m-1\right)^{2}}
\]
 on both sides of (\ref{eq:20130605-1}). Since 
\[
E^{\mathbb{Q}}\left[u_{0}^{h}\left(X_{T}\right)\right]=E^{\mathbb{Q}}\left[u_{0}^{h}\left(X_{T}\right)I_{\left\{ u_{0}\left(X_{T}\right)\geq1\right\} }\right]+E^{\mathbb{Q}}\left[u_{0}^{h}\left(X_{T}\right)I_{\left\{ u_{0}\left(X_{T}\right)<1\right\} }\right],
\]
we can apply monotone convergence theorem and bounded convergence
theorem on these two terms respectively. Therefore, (\ref{eq:20130605-1})
holds also for 
\[
h=m-\sqrt{1-\left(n-1\right)\left(m-1\right)^{2}}\mbox{ and }h=m+\sqrt{1-\left(n-1\right)\left(m-1\right)^{2}}.
\]
Hence we obtained Corollary 3.

Note that when $m=1$ in Corollary 3, $h$ is allowed to approach
to $0^+$. Since $u_{0}>0$, by Jensen's inequality, it holds that 
\begin{eqnarray*}
\lim_{h\rightarrow 0^+}\frac{E^{\mathbb{Q}}\left[u_{0}^{h}\left(X_{T}\right)\right]-u^{h}\left(T,x\right)}{Th} & = & \lim_{h\rightarrow 0^+}\frac{E^{\mathbb{Q}}\left[u_{0}^{h}\left(X_{T}\right)-1\right]}{Th}-\lim_{h\rightarrow0}\frac{u^{h}\left(T,x\right)-1}{Th}\\
 & = & \frac{E^{\mathbb{Q}}\left[\log u_{0}\left(X_{T}\right)\right]-\log u\left(T,x\right)}{T}\\
 & = & \frac{E^{\mathbb{Q}}\left[\log\frac{u_{0}\left(X_{T}\right)}{u\left(T,x\right)}\right]}{T}\\
 & \leq & \frac{1}{T}\log\frac{E^{\mathbb{Q}}\left[u_{0}\left(X_{T}\right)\right]}{u\left(T,x\right)}\\
 & \leq & \frac{1}{T}\log\frac{\left\Vert u_{0}\right\Vert _{\infty}}{u\left(T,x\right)}
\end{eqnarray*}
Therefore, by letting $h\rightarrow0$ in (\ref{eq:20130606-1-1}),
we have Corollary 4.

\subsubsection{Negative finite variation part}

Now we proceed from fundamental decompositions in a different way
in order to get a global gradient bound valid for a wider class of
$m$. From (\ref{eq:20120906-4}) and (\ref{eq:20130529-1}), we have
the global version of two fundamental decompositions. 
\begin{eqnarray}
 &  & dH_{t}\nonumber \\
 & \geq & A_{t}d\tilde{W}_{t}\nonumber \\
 &  & -\frac{1}{2}mn\left(m-1\right)^{2}\left(hf\right)^{\frac{2m-3h-2}{h}}\left|\nabla f\right|^{4}dt \\
 &  & -m\left(\left(2m-h+1-2\lambda\right)\left(m-h-1\right)+2\left(\lambda-1\right)\left(\lambda-m\right)\right)\left(hf\right)^{\frac{2m-3h-2}{h}}\left|\nabla f\right|^{4}dt.\nonumber\label{eq:global-1}
\end{eqnarray}
\begin{equation}
dY_{t}=\sqrt{2m}\left(hf\right)^{\frac{m-1}{2h}}\frac{\partial f}{\partial x^{\alpha}}d\tilde{W}_{t}^{\alpha}+m\left(2\lambda-m+h\right)\left(hf\right)^{\frac{m-h-1}{h}}\left|\nabla f\right|^{2}dt.\label{eq:global-2}
\end{equation}
 For simplicity, set 
\[
L_{1}=m\left(\left(2m-h+1-2\lambda\right)\left(m-h-1\right)+2\left(\lambda-1\right)\left(\lambda-m\right)\right)+\frac{1}{2}mn\left(m-1\right)^{2},
\]
Then (\ref{eq:global-1}) becomes 
\begin{eqnarray}
dH_{t} & \geq & A_{t}d\tilde{W}_{t}-L_{1}\left(hf\right)^{\frac{2m-3h-2}{h}}\left|\nabla f\right|^{4}dt\nonumber \\
 & = & A_{t}d\tilde{W}_{t}-L_{1}h^{-1}Y_{t}^{-1}H_{t}^{2}dt.\label{eq:20130607-3}
\end{eqnarray}
Previously, in order to obtain the submartingale property of $H$,
we always assume $L_{1}\leq0$, resulting in a constraint on $m$.
In this section, we consider the situation when 
\[
L_{1}>0.
\]
Instead of deriving a gradient bound directly, we seek for  an integral
inequality satisfied by the gradient bound as a function of $t$.
The arguments are based on the following estimate.
\begin{lem}
$H_{\cdot}e^{\int_{0}^{\cdot}L_{1}h^{-1}Y_{r}^{-1}H_{r}dr}$ is a
submartingale. Moreover,
\begin{equation}
m\left(2\lambda-m+h\right)E^{\mathbb{Q}}\left[\left.\int_{s}^{T}H_{t}\right|\mathcal{F}_{s}\right]=E^{\mathbb{Q}}\left[\left.Y_{T}\right|\mathcal{F}_{s}\right]-Y_{s}.\label{eq:2013-11-28-1}
\end{equation}

\end{lem}
Now let us define
\[
g\left(t\right)=\left\Vert Y_{t}^{-1}H_{t}\right\Vert _{\infty}.
\]

\begin{lem}
Let $u$ be a positive and bounded solution to (\ref{eq:20130306-1})
on \textup{$\left(0,\infty\right)\times\mathbb{R}^{n}$. }Assume\textup{
$L_{1}\geq0$ }and\textup{ $2\lambda-m+h>0$. Then 
\begin{equation}
g\left(s\right)\leq\frac{\left\Vert E^{\mathbb{Q}}\left[\left.\frac{Y_{T}}{Y_{s}}\right|\mathcal{F}_{s}\right]-1\right\Vert _{\infty}}{m\left(2\lambda-m+h\right)\int_{s}^{T}e^{-\int_{s}^{t}L_{1}h^{-1}g\left(u\right)du}dt}.\label{eq:20130530-1}
\end{equation}
}\end{lem}
\begin{proof}
Since $\left(e^{\int_{0}^{t}L_{1}h^{-1}Y_{r}^{-1}H_{r}dr}H_{t}\right)_{t\in\left[0,T\right]}$
is a submartingale by Lemma 13, for any $0\leq s<t\leq T$
\[
E^{\mathbb{Q}}\left[\left.e^{\int_{0}^{t}L_{1}h^{-1}Y_{r}^{-1}H_{r}dr}H_{t}\right|\mathcal{F}_{s}\right]\geq e^{\int_{0}^{s}L_{1}h^{-1}Y_{r}^{-1}H_{r}dr}H_{s},
\]
which implies
\[
E^{\mathbb{Q}}\left[\left.e^{\int_{s}^{t}L_{1}h^{-1}Y_{r}^{-1}H_{r}dr}H_{t}\right|\mathcal{F}_{s}\right]\geq H_{s}.
\]
Hence 
\[
E^{\mathbb{Q}}\left[\left.H_{t}\right|\mathcal{F}_{s}\right]\geq e^{-\int_{s}^{t}L_{1}h^{-1}g\left(r\right)dr}H_{s}.
\]
Then we have 
\[
E^{\mathbb{Q}}\left[\left.\int_{s}^{T}H_{t}dt\right|\mathcal{F}_{s}\right]\geq H_{s}\int_{s}^{T}e^{-\int_{s}^{t}L_{1}h^{-1}g\left(r\right)dr}dt.
\]
Together with (\ref{eq:2013-11-28-1}), it follows that 
\[
m\left(2\lambda-m+h\right)H_{s}\int_{s}^{T}e^{-\int_{s}^{t}L_{1}h^{-1}g\left(r\right)dr}dt\leq E^{\mathbb{Q}}\left[\left.Y_{T}\right|\mathcal{F}_{s}\right]-Y_{s}
\]
which leads to 
\[
g\left(s\right)\leq\frac{\left\Vert E^{\mathbb{Q}}\left[\left.\frac{Y_{T}}{Y_{s}}\right|\mathcal{F}_{s}\right]-1\right\Vert _{\infty}}{m\left(2\lambda-m+h\right)\int_{s}^{T}e^{-\int_{s}^{t}L_{1}h^{-1}g\left(u\right)du}dt}.
\]

\end{proof}
Let us solve this integral inequality. Set
\[
G\left(s\right)=\int_{s}^{T}e^{-\int_{s}^{t}L_{1}h^{-1}g\left(u\right)du}dt.
\]
Then 
\begin{equation}
g\left(0\right)\leq\frac{E^{\mathbb{Q}}\left[\frac{Y_{T}}{Y_{0}}\right]-1}{m\left(2\lambda-m+h\right)G\left(0\right)},\label{eq:20130610-4}
\end{equation}
 $G\left(T\right)=0$. Moreover, for any $s\in\left(0,T\right)$
\begin{eqnarray*}
G^{\prime}\left(s\right) & = & -1+L_{1}h^{-1}g\left(s\right)\int_{s}^{T}e^{-\int_{s}^{t}L_{1}h^{-1}g\left(u\right)du}dt\\
 & \leq & -1+L_{1}h^{-1}\frac{\left\Vert E^{\mathbb{Q}}\left[\left.\frac{Y_{T}}{Y_{s}}\right|\mathcal{F}_{s}\right]-1\right\Vert _{\infty}}{m\left(2\lambda-m+h\right)}.
\end{eqnarray*}
As a consequence, 
\begin{eqnarray*}
G\left(0\right) & = & G\left(T\right)-\int_{0}^{T}G\left(s\right)^{\prime}ds\\
 & \geq & G\left(T\right)+\int_{0}^{T}\frac{m\left(2\lambda-m+h\right)-L_{1}h^{-1}\left\Vert E^{\mathbb{Q}}\left[\left.\frac{Y_{T}}{Y_{s}}\right|\mathcal{F}_{s}\right]-1\right\Vert _{\infty}}{m\left(2\lambda-m+h\right)}ds\\
 & = & \frac{1}{m\left(2\lambda-m+h\right)}\int_{0}^{T}m\left(2\lambda-m+h\right)-L_{1}h^{-1}\left\Vert E^{\mathbb{Q}}\left[\left.\frac{Y_{T}}{Y_{s}}\right|\mathcal{F}_{s}\right]-1\right\Vert _{\infty}ds.
\end{eqnarray*}
This, together with (\ref{eq:20130610-4}) yields 
\[
g\left(0\right)\leq\frac{E^{\mathbb{Q}}\left[\frac{Y_{T}}{Y_{0}}\right]-1}{\int_{0}^{T}m\left(2\lambda-m+h\right)-L_{1}h^{-1}\left\Vert E^{\mathbb{Q}}\left[\left.\frac{Y_{T}}{Y_{s}}\right|\mathcal{F}_{s}\right]-1\right\Vert _{\infty}ds}
\]
as long as 
\[
\int_{0}^{T}m\left(2\lambda-m+h\right)-L_{1}h^{-1}\left\Vert E^{\mathbb{Q}}\left[\left.\frac{Y_{T}}{Y_{s}}\right|\mathcal{F}_{s}\right]-1\right\Vert _{\infty}ds>0.
\]
This means
\[
u^{m+h-3}\left|\nabla u\right|^{2}\left(T,x\right)\leq\frac{E^{\mathbb{Q}}\left[u_{0}^{h}\left(X_{T}\right)\right]-u^{h}\left(T,x\right)}{\int_{0}^{T}mh\left(2\lambda-m+h\right)-L_{1}\left\Vert E^{\mathbb{Q}}\left[\left.\frac{u_{0}^{h}\left(X_{T}\right)}{u^{h}\left(T-s,X_{s}\right)}\right|\mathcal{F}_{s}\right]-1\right\Vert _{\infty}ds}.
\]
Define 
\[
u_{\min}^{h,T}=\inf_{\left(t,x\right)\in\left[0,T\right]\times\mathbb{R}^{n}}u^{h}\left(t,x\right).
\]
It then follows that when
\begin{equation}
mh\left(2\lambda-m+h\right)-L_{1}\left(\frac{\left\Vert u_{0}^{h}\right\Vert _{\infty}}{u_{\min}^{h,T}}-1\right)>0,\label{eq:20131211-3}
\end{equation}
we have 
\begin{equation}
u^{m+h-3}\left|\nabla u\right|^{2}\left(T,x\right)\leq\frac{\left\Vert u_{0}^{h}\right\Vert _{\infty}-u^{h}\left(T,x\right)}{T\left(mh\left(2\lambda-m+h\right)-L_{1}\left(\frac{\left\Vert u_{0}^{h}\right\Vert _{\infty}}{u_{\min}^{h,T}}-1\right)\right)}.\label{eq:20131211-2}
\end{equation}
To maximize $mh\left(2\lambda-m+h\right)-L_{1}\left(\frac{\left\Vert u_{0}^{h}\right\Vert _{\infty}}{u_{\min}^{h,T}}-1\right)$,
we should take $\lambda=\frac{2m-h}{2}+\frac{h}{2}\left(\frac{\left\Vert u_{0}^{h}\right\Vert _{\infty}}{u_{\min}^{h,T}}-1\right)^{-1}$.
Then (\ref{eq:20131211-2}) becomes
\begin{eqnarray*}
 &  & u^{m+h-3}\left|\nabla u\right|^{2}\left(T,x\right)\\
 & \leq & \frac{\left\Vert u_{0}^{h}\right\Vert _{\infty}-u^{h}\left(T,x\right)}{mT\left(mh+\left(mh-m+1-\frac{h^{2}}{2}-\frac{1}{2}n\left(m-1\right)^{2}\right)\left(\frac{\left\Vert u_{0}^{h}\right\Vert _{\infty}}{u_{\min}^{h,T}}-1\right)+\frac{h^{2}}{2}\left(\frac{\left\Vert u_{0}^{h}\right\Vert _{\infty}}{u_{\min}^{h,T}}-1\right)^{-1}\right)},
\end{eqnarray*}
and the constraint (\ref{eq:20131211-3}) becomes 
\[
mh+\left(mh-m+1-\frac{h^{2}}{2}-\frac{1}{2}n\left(m-1\right)^{2}\right)\left(\frac{\left\Vert u_{0}^{h}\right\Vert _{\infty}}{u_{\min}^{h,T}}-1\right)+\frac{h^{2}}{2}\left(\frac{\left\Vert u_{0}^{h}\right\Vert _{\infty}}{u_{\min}^{h,T}}-1\right)^{-1}>0.
\]
By rewriting  them in terms of $U$, $U_{\pm}$ and $h_{\pm}$,
we have Theorem 5.

\end{document}